\documentclass[12pt]{article}
\usepackage{amsmath,amssymb,amsthm}
\usepackage[lf]{Baskervaldx}
\usepackage[bigdelims,vvarbb]{newtxmath}
\usepackage[cal=boondoxo]{mathalfa}

\usepackage{helvet}
\usepackage{microtype}
\usepackage[headings]{fullpage}
\usepackage[hidelinks]{hyperref}
\renewcommand{\thefootnote}{\fnsymbol{footnote}}
\numberwithin{equation}{section}
\DeclareRobustCommand{\PDEkeywords}[1]{\par\medskip\noindent\textbf{Keywords.} #1}
\newtheorem{theorem}{Theorem}[section]
\newtheorem{lemma}[theorem]{Lemma}
\newtheorem{proposition}[theorem]{Proposition}
\newtheorem{corollary}[theorem]{Corollary}
\theoremstyle{remark}
\newtheorem{remark}[theorem]{Remark}
\newcommand{\R}{\mathbb R}
\newcommand{\HH}{\mathcal H}
\newcommand{\hd}{\dim_{\HH}}
\newcommand{\px}{\pi_x}
\newcommand{\pt}{\pi_t}
\newcommand{\dist}{\operatorname{dist}}
\title{Generic regularity for the classical obstacle problem in dimension five}
\author{Sunghoon Kim}
\date{}
\begin{document}
\setlength{\baselineskip}{17pt}
\maketitle
\markboth{Sunghoon Kim}{Generic Regularity for the Obstacle Problem in Dimension Five}
{\renewcommand{\thefootnote}{}\footnote{Department of Mathematics, The Catholic University of Korea, 43 Jibong-ro, Bucheon-si, Gyeonggi-do 14662, Republic of Korea. E-mail: \href{mailto:math.s.kim@catholic.ac.kr}{math.s.kim@catholic.ac.kr}.}}
\begin{abstract}
Let $u^t$ be a continuous monotone family of solutions to the classical obstacle problem in a ball in $\mathbb R^5$.
We assume that, for each parameter, the boundary values at larger parameters increase at least linearly on compact portions where the boundary value at that parameter is positive.
We prove that the set of parameters for which $u^t$ has an interior singular free-boundary point has Lebesgue measure zero.
The proof uses a harmonic barrier outside a codimension-two coincidence set, comparison with a thin obstacle solution at cubic and half-integer frequencies, and a weighted spherical identity for the nonharmonic cubic.
\PDEkeywords{Obstacle problem; free boundary; singular set; monotone families; epiperimetric inequality.}
\end{abstract}


\medskip
\section{Introduction and reduction to the critical strata}\label{sec:reduction}

\subsection{The problem and the main result}
We consider the classical obstacle problem
\begin{equation}\label{eq:classical}
  \Delta u^t=\mathbf 1_{\{u^t>0\}},\qquad u^t\ge0
  \qquad\text{in } B_1\subset\mathbb R^5,
\end{equation}
where $t$ varies in $I=(-1,1)$ and the equation is understood in the sense of distributions.
The family is nondecreasing in $t$, jointly continuous in the interior, and its boundary trace is continuous in $t$ in $L^\infty(\partial B_1)$.
As in \cite[(1.2), (6.1)]{FRS}, we further assume that for each $t\in I$ and each compact set $K\subset\partial B_1\cap\{u^t>0\}$ there is $c_{t,K}>0$ such that
\begin{equation}\label{eq:boundary-increase}
 \inf_{x\in K}\bigl(u^s(x)-u^t(x)\bigr)
 \ge c_{t,K}(s-t) \qquad (t<s<1).
\end{equation}
For example, one may take Dirichlet data $g+t$ on an interval where these data are nonnegative.
Here $t$ labels an elliptic family; no evolution equation in $t$ is assumed.

Let $\Gamma^t=\partial\{u^t>0\}\cap B_1$.
Recall that a point of $\Gamma^t$ is singular if its quadratic blow-up is a nonnegative quadratic polynomial, rather than the regular half-space profile $\frac12(x\cdot e)_+^2$; see \cite[Section 3]{FRS}.
Our main result is as follows.

\begin{theorem}\label{thm:main}
Let $\{u^t\}_{t\in I}$ satisfy \eqref{eq:classical} and the monotonicity and continuity assumptions above, including \eqref{eq:boundary-increase}.
Then, for Lebesgue-almost every $t\in I$, the interior free boundary $\Gamma^t$ has no singular points.
\end{theorem}

To put this result in context, write $\Sigma^t$ for the singular set of $u^t$ and $\Sigma=\{(x,t):x\in\Sigma^t\}$.
We denote the spatial and parameter projections by $\px$ and $\pt$.
At a point $(x,t)\in\Sigma$, the quadratic blow-up has the form $p_{2,x,t}(h)=\frac12 h\cdot M_{x,t}h$, where $M_{x,t}\ge0$ and $\operatorname{tr}M_{x,t}=1$.
Let $\Sigma_m$ denote the set of such points with $\dim\ker M_{x,t}=m$.

Monneau \cite{Monneau} proved the planar case of Schaeffer's conjecture \cite{Schaeffer} for increasing Dirichlet data.
Figalli, Ros-Oton, and Serra \cite[Theorems 1.1 and 9.5]{FRS} showed that the singular set is empty for almost every parameter in dimensions at most four.
In dimension five their estimate gives $\HH^1(\Sigma^t)=0$ for almost every $t$.
Theorem~\ref{thm:main} concerns the different set $\pt\Sigma$: it states that even a single singular point can occur only at parameters in a Lebesgue null set.

Two further results enter our reduction.
Franceschini and Zato\'n \cite[Theorems 1.2 and 7.9]{FZ} separate the set of singular pairs with expansions of every order, denoted by $\Sigma^\infty$, from a set with smaller spatial projection.
Their estimates give
\begin{equation*}
 \hd(\pt\Sigma^\infty)=0,\qquad
 \hd\bigl(\px(\Sigma\setminus\Sigma^\infty)\bigr)\le n-2.
\end{equation*}
The second estimate, which concerns $x$ alone, is useful here only where we can compare the parameters of nearby singular pairs.
For homogeneous Signorini solutions, Franceschini and Savin \cite[Theorem 1]{FrSa} exclude frequencies in $(2k,2k+1)$.
This also applies to the nonintegral classical blow-up limits used below, since their part odd in the normal variable is harmonic and must vanish at nonintegral homogeneity.

For comparison, Fern\'andez-Real and Torres-Latorre \cite{FRTL} prove generic regularity for the Signorini problem in ambient dimensions up to four; their treatment of cubic points uses the thin-obstacle convergence theorem of Savin and Yu \cite{SYinteger}.
Carducci and Colombo \cite{CC} consider the fractional obstacle problem for $0<s<1$.
These results concern families that already satisfy a thin or fractional obstacle equation.
Here the equation is \eqref{eq:classical}; an exact Signorini solution is used to approximate the rescaled remainder, and the approximation error must be controlled in the original equation.

We use the stratification and comparison framework of \cite[Sections 6--9]{FRS}.
Four cases require estimates at a limiting exponent.
The frequency-four case without a completed quartic expansion is treated in Proposition~\ref{prop:failed-quartic}.
For the ordinary stratum with a three-dimensional quadratic kernel, Section~\ref{sec:ordinary} constructs a barrier outside a narrow tube.
The parameter projection is then controlled by a logarithmic estimate covered by Altaf and Cs\H{o}rnyei \cite[Theorem 2]{AltafCs}.
For the nonharmonic cubic, Section~\ref{sec:cubic} first improves a thin-obstacle approximation in the even part and then estimates the parameter increment by a weighted integral on the sphere.
For the harmonic cubic with next frequency $7/2$, Section~\ref{sec:slit} uses three independent directions of neighboring singular points to improve the approximation by a fixed slit.
Once the slit is determined, a barrier in the positivity set gives the required parameter estimate.
The proof of Theorem~\ref{thm:main} is completed in Section~\ref{sec:conclusion}.

\subsection{Dimension bounds and estimates for the parameter projection}\label{sec:projection-reduction}
For $E\subset\Sigma$, forward cleaning of order $s>0$ means the following: for every $(x_0,t_0)\in E$ and every $\varepsilon\in(0,s)$ there exists $\rho=\rho(x_0,t_0,\varepsilon)>0$ such that
\begin{equation*}
 (x,t)\in E,\quad |x-x_0|<\rho,\quad 0<t-t_0<\rho
 \quad\Longrightarrow\quad
 t-t_0\le |x-x_0|^{s-\varepsilon}.
\end{equation*}
Constants can be absorbed by decreasing the exponent.
If $\hd(\px E)\le\beta<s$, then $\hd(\pt E)\le\beta/s$ by \cite[Corollary 7.8(a)]{FRS}.
At equality $\beta=s$, one can instead use a pointwise modulus $\omega(r)r^s$, $\omega(r)\to0$, and $\HH^s$-$\sigma$-finiteness of $\px E$ to obtain $\HH^1(\pt E)=0$ by \cite[Proposition 7.7(a)]{FRS}.
The same principles apply backward after reversing the parameter.
Equality of the dimensions alone does not give the conclusion.

We record explicitly the comparison used throughout the argument.
Let $P$ be a polynomial with $\Delta P=1$ and let $e(r)\ge\sup_{B_r}(P-u^{t_0})$.
In coordinates with transverse axis $x_5$, set
\begin{equation*}
\begin{aligned}
 \psi_r(x)&=\tfrac12+\frac{5x_5^2-|x|^2}{r^2},\\
 D_r&=\{x\in\partial B_r:\ |x_5|>r/\sqrt{10}\},\\
 m_r&=\inf_{D_r}\bigl(u^t-u^{t_0}\bigr),\\
 b(x)&=P(x)+\frac{2}{9}m_r\psi_r(x)-e(r).
\end{aligned}
\end{equation*}
On $\partial B_r\setminus D_r$ the term involving $\psi_r$ is nonpositive.
On $D_r$ it is at most $m_r$ since $\max_{\partial B_r}\psi_r=9/2$.
Consequently $b\le u^t$ on $\partial B_r$.
Since $\Delta\psi_r=0$ and $\Delta u^t\le1$, the difference $u^t-b$ is superharmonic and nonnegative in $B_r$.
In particular,
\begin{equation}\label{eq:comparison}
 u^t\ge P+\frac{m_r}{18}-e(r)\qquad\text{in } B_{r/2},
\end{equation}
where $\psi_r\ge1/4$ in $B_{r/2}$.
If $P$ has a lower bound only on the old coincidence set, it suffices to apply \eqref{eq:comparison} there; at all other points $u^t>0$ by monotonicity.
This is the signed harmonic comparison of \cite[Lemma 9.1]{FRS}, with the actual error $e(r)$ retained.

\paragraph{Lower strata.}
For $m\le2$ the spatial bounds in \cite[Proposition 8.1]{FRS} are $m-1$ on the anomalous set (the lowest stratum is treated separately there) and $m$ on its complement.
The respective parameter exponents in \cite[Lemma 9.3(a)--(b)]{FRS} can be taken as $2-\varepsilon$ and $3-\varepsilon$, which are strictly larger than the corresponding spatial dimensions for sufficiently small $\varepsilon>0$.
At $m=3$, the anomalous limit has restriction $-\sum_{i=1}^3\nu_i h_i^2$ to the spine, where $\nu>0$, $|\nu_i|\le\nu$ and $\sum_i\nu_i=2\nu$ by \cite[Proposition 3.9(a), (3.11)]{FRS}.
All $\nu_i$ are nonnegative: otherwise the sum would be less than $2\nu$.
At least two are positive, so the zero set on the spine has dimension at most one.
Apply the density extraction of \cite[Lemma 7.6]{FRS} as in \cite[Proposition 8.1(a)]{FRS}: a spatial dimension greater than one would produce, for some $\beta>1$, an accumulation set of positive $\HH^\beta$-measure for a single anomalous blow-up.
\cite[Lemma 6.4]{FRS} places it in that zero set, a contradiction.
Thus the parameter dimension of the anomalous $m=3$ stratum is at most $1/2$.

For ordinary points of $\Sigma_3$, the second frequency is at least three, and is an integer by \cite[Proposition 3.9(a)]{FRS}.
When it is at least four, \cite[Lemmas 3.6--3.7]{FRS} give $\|u-p_2\|_{L^\infty(B_r)}\le Cr^4$.
For any $\varepsilon\in(0,1)$, the lower bound of \cite[Lemma 9.2(a)]{FRS} on the transverse caps is $c_\varepsilon r^\varepsilon(t-t_0)$.
Inserting it into \eqref{eq:comparison} with $P=p_2\ge0$ gives the exponent $4-\varepsilon>3$.
The argument in Section~\ref{sec:ordinary} treats the remaining case $\lambda_2=3$ and uses only the upper remainder bound $Cr^3$, so it also applies if the solution is exactly quadratic in a neighborhood.

\paragraph{Exclusion of the top subcubic stratum.}\label{par:subcubic}
Suppose a point belongs to $\Sigma_4^{<3}$.
By \cite[Proposition 3.9(b)]{FRS}, a normalized second blow-up is a nonzero global homogeneous Signorini solution $q$, with homogeneity $2+\alpha_\circ\le\lambda<3$ for some $\alpha_\circ>0$.
Its part odd in the normal variable is an entire harmonic function: the nonpositive measure $\Delta q$ is supported on the thin plane and cancels under odd reflection.
An entire harmonic homogeneous function of noninteger degree is zero, by its harmonic Taylor expansion at the origin.
Thus $q$ is even and \cite[Theorem 1]{FrSa} excludes its homogeneity in $(2,3)$.
Consequently $\Sigma_4^{<3}=\emptyset$.

\paragraph{Top nonharmonic cubic.}
At a point in $\Sigma_4^{\ge3}\setminus\Sigma_4^{3\mathrm{rd}}$, the unique cubic limit has even part
\begin{equation*}
 -|z|\left(x'\cdot A x'-\tfrac13(\operatorname{tr}A)z^2\right),
 \qquad A\ge0,\quad A\ne0,
\end{equation*}
after choosing the normalization of $A$; see \cite[Lemma 5.2, Proposition 5.4 and (6.27)]{FRS}.
The proof of \cite[Lemma 6.6]{FRS} actually places every limiting spatial secant between points of this nonharmonic cubic set in the kernel of $A$ at the base point: the translation invariance of the even part in \cite[Lemma 6.5]{FRS} gives $Av=0$.
For rank$(A)\ge2$, these kernels have dimension at most two.
Repeat \cite[Proposition 8.3]{FRS} with the full kernel padded to a two-plane: localize the admissible radius and approximate the tangent plane by one of countably many fixed planes.
On each piece the inverse of orthogonal projection to the fixed plane is Lipschitz.
The projection is therefore covered by countably many Lipschitz two-graphs, hence $\HH^2$-$\sigma$-finite.
The backward little-$o(r^2)$ estimate of \cite[Lemma 9.4]{FRS} and \cite[Proposition 7.7(a)]{FRS} then give null parameter projection.
Only rank$(A)=1$ remains in this cubic stratum.

\paragraph{Top harmonic cubic jet and noninteger residual.}
Let $\mathcal P$ be the fourth-order almost-square ansatz of \cite[Definition 4.5]{FRS}.
In normal coordinates it has the form $\mathcal P=\tfrac12(z+a)^2+zQ_3$, where $a=p_3/z$ is homogeneous quadratic and $Q_3$ is homogeneous cubic.
Completing the square once more gives the global bound
\begin{equation}\label{eq:ansatz-lower}
 \mathcal P=\tfrac12(z+a+Q_3)^2-aQ_3-\tfrac12Q_3^2
 \ge -Cr^5\qquad\hbox{in }B_r,\quad 0<r<1.
\end{equation}
This estimate does not require restriction to the coincidence set.
\cite[Proposition 5.8]{FRS} gives residual frequency strictly above three after the harmonic cubic.
On $\Sigma_4^{3\mathrm{rd}}\setminus\Sigma_4^{\ge4}$ the residual frequency $\lambda_3$ lies in $(3,4)$ by \cite[Proposition 4.12(a)]{FRS}.
In the proof of \cite[Lemma 8.4(a) and Remark 8.5]{FRS}, three independent limiting directions and local frequency matching force, by \cite[Lemma 6.8]{FRS}, a two-dimensional Signorini profile.
The two-dimensional list in \cite[Lemma A.2]{FRS} permits only $7/2$ in $(3,4)$.
To apply \cite[Proposition 7.3]{FRS}, fix a base point with frequency $\lambda_3\in(3,4)\setminus\{7/2\}$.
If the required two-plane approximation failed for a fixed aperture, one could choose radii $r_j\downarrow0$ and three points at each radius whose frequencies tend to that of the base point.
Choose the three secants successively outside the aperture tubes about the spans of the preceding ones.
After rescaling by $r_j$, the limiting secants are linearly independent.
A common subsequence of the base blow-ups then has these three translation directions, contrary to the two-dimensional classification at $\lambda_3$.
This supplies the frequency-localized two-plane approximation required in \cite[Proposition 7.3]{FRS}.
Thus the spatial dimension is at most $n-3=2$.
For any $3<\gamma<\lambda_3$, use $\|u-\mathcal P\|_\infty=O(r^\gamma)$ by \cite[Lemmas 4.7 and 4.9]{FRS}, $\mathcal P\ge-Cr^5$, and the $crt$ perturbation lower bound to obtain order $\gamma-1>2$, as in \cite[Lemma 9.3(d)]{FRS}.
Countable localization of $\gamma$ suffices for the slicing argument.
Only the exact $7/2$ case remains in this interval; its known order $5/2$ is smaller than the possible spatial dimension three.

\paragraph{Completed quartic jet.}
On $\Sigma_4^{4\mathrm{th}}$, \cite[Proposition 4.15]{FRS} excludes residual frequency four, and \cite[Proposition 4.12(b)]{FRS} together with \cite[Theorem 1]{FrSa} excludes $(4,5)$.
Thus the completed quartic points are in $\Sigma^*$ of \cite{FRS} and have order-$4-\varepsilon$ cleaning for every $\varepsilon>0$ by \cite[Lemma 9.3(e)]{FRS}.
The infinite-expansion portion has null parameter projection by \cite[Theorem 1.2]{FZ}.
Its complement has spatial dimension at most three by \cite[Theorem 7.9]{FZ}, hence parameter dimension at most $3/4$ by sending $\varepsilon$ to zero.
The failed quartic jet is a different stratum, treated next.

\subsection{The frequency-four case without a completed quartic expansion}\label{sec:failed-quartic}
\begin{proposition}\label{prop:failed-quartic}
Let $E=\Sigma_4^{\ge4}\setminus\Sigma_4^{4\mathrm{th}}$ in the notation of \cite{FRS}.
For the families above, $\HH^1(\pt E)=0$.
\end{proposition}
\begin{proof}
Translate a point of $E$ to $(0,0)$, let $L=\{z=0\}$ with $x=(x',z)$, and write $w=u^0-\mathcal P$ and $u_r(x)=u^0(rx)$.
\cite[(6.45), proof of Lemma 6.12]{FRS} gives local uniform subsequential compactness of $r^{-4}w(r\cdot)$, with each accumulation point a harmonic homogeneous quartic of nonnegative trace on $L$ and nonzero even part.
An even harmonic quartic is determined by its trace: successive coefficients of $z^{2j}$ follow from $\partial_{zz}=-\Delta'$.
Therefore the trace is nonzero.

The moment computation in the proof of \cite[Lemma 4.14]{FRS} can be applied to $w$ without the completed quartic jet.
Indeed $\Delta[w(r\cdot)]=-r^2\mathbf 1_{\{u_r=0\}}$, and \cite[(6.45)]{FRS} at scale $2r$ gives $\|w(r\cdot)\|_{L^\infty(B_2)}=O(r^4)$.
Test this distributional equation with a nonnegative smooth cutoff supported in $B_2$ and equal to one on $B_1$.
Integration by parts gives $|\{u_r=0\}\cap B_1|\le Cr^2$.
Moreover, on this set, the equality
\begin{equation*}
 0=u_r(x)=\tfrac12 r^2z^2+r^3 p_3(x)+O(r^4)
\end{equation*}
and the divisibility of $p_3$ by $z$ imply $|z|\le Cr$.
For a homogeneous harmonic quartic $H$ with $H|_L\ge0$, Green's identity gives
\begin{equation*}
 \frac{d}{dr}\left(r^{-4}\int_{\partial B_1}w(rx)H(x)\,dS_x\right)
  =-r^{-3}\int_{B_1\cap\{u_r=0\}}H(x)\,dx\le C_H,
\end{equation*}
where $H(x',z)\ge-C_H|z|$ was used.
The moment is bounded, and after subtracting $C_Hr$ is nonincreasing in $r$.
Thus it has a finite limit at zero.
Such tests span all harmonic quartics: add a sufficiently large multiple of the harmonic even extension of $f(x')=|x'|^4$, namely $f-z^2\Delta'f/2+z^4(\Delta')^2f/24$, to a given test to make its trace nonnegative.
It follows that the limits of the moments determine a unique accumulation point $Q$, by testing the difference of two candidate limit quartics against itself.
In particular
\begin{equation*}
 \|u^0-\mathcal P-Q\|_{L^\infty(B_r)}=o(r^4).
\end{equation*}

On $\{u^0=0\}\cap B_r$ one has $|z|\le Cr^2$ in physical coordinates, so $Q\ge -Cr^5$ there by homogeneity and its nonnegative flat trace.
The same lower bound for $\mathcal P$ follows from the almost-square formula.
Because $\Delta(\mathcal P+Q)=1$, apply \cite[Lemmas 9.1 and 9.2(b)]{FRS} to $\mathcal P+Q$, with the actual scale error $o(r^4)$ as the comparison constant in the proof of \cite[Lemma 9.1]{FRS}.
On $B_{r/2}\cap\{u^0=0\}$ this gives $u^t\ge crt-o(r^4)$ for $t>0$.
Elsewhere positivity is inherited by monotonicity.
Taking $r=3|x|$ yields forward cleaning by a point-dependent modulus $\omega(|x|)|x|^3$, with $\omega(r)\to0$.

Choose a unit vector $v\in L$ with $Q(v)>0$.
By monotonicity, all later singular spatial points belong to the old contact set.
Their limiting directions lie in $L$ and cannot approach $v$ or $-v$, because $\mathcal P\ge-Cr^5$ and $Q(v)=Q(-v)>0$.
By compactness this excludes a symmetric cone about the two-plane $\Pi=\operatorname{span}\{v,e_z\}$ for all such secants from the base point at a sufficiently small radius.
Partition $E$ countably by lower bounds on this radius and aperture, local spatial diameter smaller than the radius bound, and a countable choice of approximate $\Pi$ in the Grassmannian (narrow the common cone).
For any pair in one piece take the earlier point as base.
The symmetric exclusion therefore holds for the secant in either orientation, and for equal times as well by the same contact implication.
Consequently projection to the fixed $\Pi^\perp$ has a Lipschitz inverse on the spatial projection of each piece.
The pieces lie in countably many Lipschitz three-graphs, and their projections are $\HH^3$-$\sigma$-finite.
Apply \cite[Proposition 7.7(a)]{FRS} with $s=\beta=3$ and the forward little-$o(r^3)$ estimate on finite-measure pieces.
\end{proof}

\section{Logarithmic projection and the ordinary rank-two stratum}\label{sec:ordinary}

\subsection{A logarithmic endpoint projection lemma}

\begin{theorem}[Special case of {\cite[Theorem 2]{AltafCs}}]\label{sard:th:sard}
Let $d\ge2$, $E\subset\R^d$, and $T:E\to\R$.
Suppose there are constants $C>0$ and $0<r_0<1$ such that for $0<|x-y|<r_0$,
\begin{equation}\label{sard:eq:modulus}
 |T(x)-T(y)|\le C |x-y|^d
    \bigl(\log(1/|x-y|)\bigr)^{d-1}.
\end{equation}
Then the Lebesgue outer measure of $T(E)$ is zero.
\end{theorem}

To see the connection with \cite[Theorem 2]{AltafCs}, consider more generally $\omega_q(r)=Cr^d(\log(1/r))^q$, for $0\le q\le d$ and small $r$.
This function is increasing and convex near zero: its first and second derivatives have leading terms $Cd r^{d-1}(\log(1/r))^q$ and $Cd(d-1)r^{d-2}(\log(1/r))^q$, respectively.
Its ratio to $r^d$ is nonincreasing.
Choose $a>0$ in this range and extend $\omega_q$ affinely for $r\ge a$ with slope $m=\omega_q'(a)>0$.
For $r\ge a$ the sign of the derivative of $\omega_q(r)/r^d$ is the sign of $rm-d[\omega_q(a)+m(r-a)]$.
This expression is nonpositive at $a$ and decreases thereafter.
The extended modulus therefore meets all the monotonicity and convexity assumptions of their theorem.
The change of variables $t=\log(1/r)$ gives
\begin{equation*}
 \int_0^a \omega_q(r)^{-1/d}\,dr
 =C^{-1/d}\int_{\log(1/a)}^\infty t^{-q/d}\,dt=\infty .
\end{equation*}
Apply the theorem on small cubes of diameter less than $a$ and less than the radius in the local modulus.
In particular it gives the asserted case $q=d-1$.
For completeness, we also include a direct proof of that case.

\begin{proof}
By countably many small cubes, it suffices to assume $E$ is bounded and $\operatorname{diam}E<r_0/2$, with $r_0$ also small enough that the modulus on the right of \eqref{sard:eq:modulus} is increasing.
On each such piece $T$ extends continuously to its compact closure, obeying the same modulus.
Thus we can assume $E$ compact without any measurability assumption on the original image.

Suppose $F=T(E)$ has positive Lebesgue measure $M_0$.
Fix $1<p<d$.
At stage $j\ge1$ remove the central mass fraction
\begin{equation*}
 \delta_j=\frac{c_0}{(j+2)(\log(j+2))^p}
\end{equation*}
from each of the stage $j-1$ intervals, splitting the remainder into a left and right child of equal mass.
Here the mass is Lebesgue measure restricted to $F$, and $c_0>0$ is fixed small enough that all $\delta_j<1/2$.
This is done using the continuous distribution function of that restricted Lebesgue measure.
The two closed child intervals are disjoint: the intervening gap contains mass $\delta_j M_{j-1}$ and hence has length at least $\delta_j M_{j-1}$.
The mass of every level-$j$ child equals
\begin{equation*}
 M_j=M_0 2^{-j}\prod_{i=1}^j(1-\delta_i)\asymp 2^{-j},
\end{equation*}
since $\sum_j\delta_j<\infty$.
All inverse images $E_{j,k}$ of the resulting closed intervals are nonempty compact sets, nested down the tree.

It follows from \eqref{sard:eq:modulus} that any two sibling sets at level $j$ have positive spatial distance at least
\begin{equation*}
 d_j=c_1 2^{-j/d}\delta_j^{1/d}(j+1)^{-(d-1)/d}
\end{equation*}
for a fixed $c_1>0$ sufficiently small.
Indeed $\log(1/d_j)\le C_d(j+1)+\log(1/c_1)+C_{d,p,c_0}$, while the time separation is at least $c2^{-j}\delta_j$.
Substituting $d_j$ in the modulus and decreasing $c_1$ makes this modulus smaller than the latter separation for every $j$.
The smallness also ensures $d_j<r_0$.
Hence distances less than $d_j$ cannot occur.

Set $r_j=d_j/4$ for $j\ge1$, and fix $r_0'>r_1$ for the root radius.
For a compact set $K$ let $K^{r}=K+B_r$ be its open $r$-neighborhood.
The level $j+1$ sibling neighborhoods of radius $r_{j+1}$ are disjoint.
Their union, thickened by $r_j-r_{j+1}$, is contained in their parent's $r_j$-neighborhood.
The Brunn--Minkowski inequality gives for $j\ge1$, with $\omega_d=|B_1|$,
\begin{equation*}
 \min_k\left(\frac{|E_{j,k}^{r_j}|}{\omega_d}\right)^{1/d}
 \ge 2^{1/d}\min_k\left(\frac{|E_{j+1,k}^{r_{j+1}}|}{\omega_d}\right)^{1/d}
        +(r_j-r_{j+1}).
\end{equation*}
The same inequality holds at the root with radius $r_0'$.
Apply it down to any finite level and omit the nonnegative last contribution.
The fixed root volume to power $1/d$ is at least a fixed positive multiple of
\begin{equation*}
 \sum_{j=1}^{N}2^{j/d}(r_j-r_{j+1})
 \gtrsim\sum_{j=1}^{N}\delta_j^{1/d}(j+1)^{-(d-1)/d}.
\end{equation*}
Here $r_{j+1}/r_j\le2^{-1/d}$.
The last sum diverges as $N\to\infty$, being comparable for large $j$ to $\sum_j j^{-1}(\log j)^{-p/d}$ with $p/d<1$.
This contradicts boundedness of the root neighborhood.
\end{proof}

\begin{remark}
The hypothesis $d\ge2$ matters: for $d=1$ the analogous exponent $d-1=0$ is only a Lipschitz condition and the identity on an interval is a counterexample.
The statement also applies under a smaller nonnegative logarithmic exponent, by monotonicity of the power of the logarithm for small radii.
\end{remark}

\subsection{A lower bound for the increment outside a codimension-two tube}

\begin{lemma}\label{sard:lem:barrier}
Let $n\ge3$, $x=(y,z)\in\R^{n-2}\times\R^2$, and consider an old coincidence set $K_0$ with
\begin{equation}\label{sard:eq:tube}
 K_0\cap B_R\subset\{|z|\le C_0 |x|^{1+\eta}\},\qquad \eta>0.
\end{equation}
Consider a family of classical obstacle solutions with the uniform boundary monotonicity of Figalli--Ros-Oton--Serra \cite[(1.2)]{FRS}.
For a fixed cone positively separated from $\{z=0\}$ and $R$ sufficiently small, its actual forward increment obeys
\begin{equation}\label{sard:eq:gain}
 u^h(x)-u^0(x)\ge c h (\log(1/|x|))^{-1/\eta}
 \qquad (h>0)
\end{equation}
on that cone near zero.
The constant can depend on the center, the solution and the cone.
\end{lemma}

\begin{proof}
Put $m=n-2$, $r=|x|$, $s=\log(1/r)$, and choose $C_1>\log C_0$.
Write
\begin{equation*}
 \tau=s-C_1/\eta,\qquad \theta=\log(r/|z|),\qquad
 v=\frac{\theta}{\eta\tau},\qquad \mathfrak{a}=\eta^{-1}.
\end{equation*}
For large $\tau$ the closed angular section $0\le v\le1$ lies in the old positive phase for $0<r<R$, because $v\le1$ gives $|z|\ge e^{C_1}r^{1+\eta}>C_0r^{1+\eta}$.
Define on $0\le v<1$ the function
\begin{align*}
 \mathcal B&=G(\tau)f(\tau,v),\qquad f(\tau,v)=(1-v)(1-v/\tau),\\
 G(\tau)&=\tau^{-\mathfrak{a}}\exp\!\left(\frac{K(\log\tau+1)}{\tau}\right).
\end{align*}
Here $K>0$ will be chosen, fixed thereafter.
The function $\mathcal B$ is positive in the interior for $\tau>1$, vanishes at $v=1$, and tends uniformly to zero at the tip.

We check the sign directly.
In polar variables $r^2\Delta=\partial_{ss}-m\partial_s+\Delta_{\mathbb S^{n-1}}$.
Since $\log|z|$ is harmonic away from $z=0$ and $\log|z|=\log r-\theta$,
\begin{equation*}
 \Delta_{\mathbb S^{n-1}}\theta=m,\qquad
 |\nabla_{\mathbb S^{n-1}}\theta|^2=e^{2\theta}-1 .
\end{equation*}
For the latter identity use $b=|z|/r$ on the sphere.
Then $|\nabla_{\mathbb S^{n-1}}b|^2=1-b^2$.
At fixed angle $v_\tau=-v/\tau$.
Denote the effective exponent by
\begin{equation*}
 \mathfrak{a}_{\mathrm{eff}}=-\tau G'/G=\mathfrak{a}+K\log\tau/\tau .
\end{equation*}
The derivative $f_\tau$ holds $v$ fixed, and $\mathcal B_{\tau\tau}$ holds the angle fixed.
Thus
\begin{equation}\label{sard:eq:sign}
 \frac{r^2\Delta \mathcal B}{G/\tau}
 =m\left[\mathfrak{a}_{\mathrm{eff}} f+(v+\mathfrak{a})f_v-\tau f_\tau\right]
   +\frac{\tau \mathcal B_{\tau\tau}}G
   +\frac{e^{2\theta}-1}{\eta^2\tau}f_{vv}.
\end{equation}
Throughout $0\le v\le1$,
\begin{align*}
 f_v&=-1-\tau^{-1}+2v/\tau,\qquad f_{vv}=2/\tau,\\
 f_\tau&=v(1-v)/\tau^2,\qquad |\tau \mathcal B_{\tau\tau}/G|\le C_{K,\eta}/\tau .
\end{align*}
The last bound follows from $G'/G=O(\tau^{-1})$, $G''/G=O(\tau^{-2})$, and the first two total radial derivatives of $f(\tau,v(\tau))$ being $O(\tau^{-1})$ and $O(\tau^{-2})$ uniformly.

For $\theta\le 2\log\tau$, one has $v\le 2\log\tau/(\eta\tau)$, so $f\ge1/2$ eventually and
\begin{equation*}
 \mathfrak{a} f+(v+\mathfrak{a})f_v-\tau f_\tau\ge -C_\eta\frac{\log\tau}{\tau}.
\end{equation*}
Choose $K>2(C_\eta+1)$ and then take $\tau$ large.
The correction in $\mathfrak{a}_{\mathrm{eff}}$ absorbs this loss as well as the $O(1/\tau)$ derivative term.
For $\theta\ge2\log\tau$, the first bracket in \eqref{sard:eq:sign} is uniformly bounded, while the last term is $2(e^{2\theta}-1)/(\eta^2\tau^2)$ and is at least of order $\tau^2$.
Hence $\Delta \mathcal B\ge0$ on the full truncated region for sufficiently small $R$.

The actual difference $w_h=u^h-u^0$ is nonnegative and harmonic on the old positive phase.
On the compact section $v\le1$ of the outer sphere $r=R$ it is bounded below by $c_0h$.
To see this, every old positive component reaches a boundary patch on which the old boundary value is positive: otherwise $\Delta u^0=1$ in that component and the maximum principle would make it nonpositive.
The uniform boundary increment \cite[(1.2)]{FRS} and a local boundary barrier or Harnack chain give a positive lower multiple of $h$ at its chosen interior point.
Harnack chains within the fixed old positive phase, and a finite subcover of the compact section, give the claimed bound on the section.
These are the same comparisons used in \cite[Lemma 9.2(a)]{FRS}, here only on a fixed radius, not on the contracting region.

Use a small multiple of $h\mathcal B$ as a lower comparison to $w_h$, first on a truncated annulus.
On the outer part it is dominated as just proved; on the lateral part it is zero; on the inner part subtract its maximum there, which tends to zero for fixed $h$ as the inner radius tends to zero.
The maximum principle applies because the difference is superharmonic.
Finally on a cone separated from $\{z=0\}$ the variable $\theta$ is bounded, hence $f\to1$ there and $G\asymp (\log(1/r))^{-\mathfrak{a}}$.
This proves \eqref{sard:eq:gain}.
\end{proof}

\subsection{Ordinary rank-two singular times in dimension five}

\begin{corollary}\label{sard:cor:closure}
In dimension five, under the family hypotheses of Theorem~\ref{thm:main}, the set of parameters carrying an ordinary lower singular point with quadratic kernel dimension three has Lebesgue measure zero.
Here ordinary means the nonanomalous part of the lower stratum, equivalently second frequency at least three, in the notation of \cite{FRS}.
\end{corollary}

\begin{proof}
At such a point center space and parameter at zero.
The positive definite normal quadratic blowup is $p_2(y,z)=z^{\mathsf T}Az/2$.
By \cite[Lemmas 3.6--3.7]{FRS}, $|u^0-p_2|\le C r^3$ on $B_r$: the lower bound three on the limiting frequency gives the exact $r^6$ upper bound for the squared spherical norm, and \cite[Lemma 3.7]{FRS} gives the local uniform estimate.
If $u^0$ agrees with its quadratic blowup locally, the same estimate holds directly.
On the old contact set this yields $|z|\le Cr^{3/2}$, i.e., \eqref{sard:eq:tube} with $\eta=1/2$.
Choose a unit normal direction in the normal two-plane.
The caps $D_r$ of \cite[Lemma 9.1]{FRS} for this direction lie in the cone of Lemma~\ref{sard:lem:barrier}; thus
\begin{equation*}
 \inf_{D_r}(u^h-u^0)\ge c h (\log(1/r))^{-2}.
\end{equation*}
Apply the harmonic polynomial comparison of \cite[Lemma 9.1]{FRS} with $p=p_2$, which has $\Delta p_2=1$ and is nonnegative.
Its harmonic cap barrier is bounded below by a positive constant on $B_{r/2}$; the subtracted error is $Cr^3$.
It follows that $u^h>0$ throughout $B_{r/2}$ whenever
\begin{equation*}
 h>C' r^3(\log(1/r))^2 .
\end{equation*}
Using $r=3|x-x_0|$ for a neighboring singular event yields a one-sided time modulus at every base point $(x_0,t_0)$ of the stratum.
With fixed spatial coordinate, $r\downarrow0$ shows that distinct times cannot both have events of this stratum at that coordinate.

The spatial projection of the ordinary stratum is locally contained in a $C^{1,1}$ three-dimensional manifold by \cite[Proposition 8.1(c)]{FRS}.
On compact subcharts we use coordinates with a bi-Lipschitz parametrization, which changes the modulus only by a constant at small distances.
Local constants and radii can be made uniform on countably many pieces: restrict to points for which the one-sided constant is at most $j$ and the admissible radius is at least $1/j$, and split these further into small chart cubes.
For any pair on such a piece, use the event with earlier parameter as base point to obtain the symmetric modulus \eqref{sard:eq:modulus} with $d=3$.
The permitted constants and radii therefore cover the entire spatial projection, without any uniform a priori estimate across singular points.
Theorem~\ref{sard:th:sard} gives a null image for each piece, and then for their countable union.
\end{proof}

\begin{remark}
This argument addresses only the displayed ordinary lower stratum.
In particular it gives no proof for any rank-one top stratum or for lower anomalous frequency-two points, for which the polynomial rate \eqref{sard:eq:tube} used here does not follow from the same input.
The proof uses only quadratic approximation with error $O(r^3)$; no quantitative rate of convergence to a cubic polynomial is required.
\end{remark}

\section{The nonharmonic cubic with rank-one even part}\label{sec:cubic}

\subsection{The cubic blow-up and its translation directions}
Consider a continuous monotone classical obstacle family of \cite[(6.1)]{FRS} in $\R^5$, with uniform increase on positive boundary portions as in \cite[(1.2)]{FRS}.
Our cohort is $\Sigma_4^{\ge3}\setminus\Sigma_4^{\mathrm{3rd}}$ with $\operatorname{rank} A=1$ in the even cubic below, not the harmonic-cubic next-frequency stratum.
At a fixed such mark, translate it to $(0,0)$ and write $x=(y,z)$ with $p_2=z^2/2$, $L=\{z=0\}$, $w=u^0-p_2$ and $K=\{u^0=0\}$.
Then
\begin{equation}\label{cubic:eq:q}
q=q_e+q_o,\qquad
q_e=|z|({\rm tr}(A)z^2/3-y^{\mathsf T}Ay),
\qquad
q_o=z({\rm tr}(B)z^2/3-y^{\mathsf T}By)
\end{equation}
is the unique, nonzero limit of $r^{-3}w(r\,\cdot)$.
Here $A\ge0$ is nonzero and $B$ is symmetric.
This uses \cite[Proposition 5.4 and Lemma 5.2]{FRS}; the latter is the even cubic classification, while the odd part is an ordinary harmonic polynomial.
\cite[Proposition 5.4]{FRS} supplies no quantitative rate.
Its weighted trace monotonicity input is \cite[Lemma 5.3]{FRS}, not a full Dirichlet epiperimetric inequality.

The exact matching premise in \cite[Lemma 6.5]{FRS} is $|x_j|\le r_j\downarrow0$ with actual neighboring rank-one quadratic marks, a normalized limit of the old residual on that same sequence, and $\lambda_j^{\mathrm{2nd}}\to\lambda^{\mathrm{2nd}}$.
Here the frequencies equal three, and raw convergence and normalized Lipschitz bounds imply the required strong interior $W^{1,2}$ convergence: gradients converge harmonically off the vanishing contact slab, and the omitted slab is estimated by the common gradient bound.
\cite[Lemma 6.6]{FRS} uses this to constrain the spatial secants to a fixed $(n-2)$-plane at each base.
In the rank-one even-cubic case this plane is $\ker A\times\{0\}$, and projected points are covered by countably many three-dimensional Lipschitz graphs.
The graph construction is the pointwise cone estimate, a countable radius localization, and a finite Grassmannian net, as in the proof of Proposition~\ref{prop:failed-quartic}.

The backward exclusion of \cite[Lemma 9.4]{FRS} bounds the parameter separation by $o(r^2)$.
Its comparison uses the qualitative cubic expansion and places a separating surface in the coincidence set.
The estimate suffices when the spatial projection is carried by two-dimensional graphs.
For a rank-one cubic the graphs are three-dimensional.
We therefore estimate the separation directly, using the sign of the difference between the earlier and later solutions.

For the genuine thin problem, Savin--Yu \cite[Theorem 1.1]{SYinteger} prove power convergence at odd integer contact frequencies, and Carducci--Velichkov \cite[Theorem 1.1]{CV} prove the odd epiperimetric inequality.
The latter applies to boundary data, not only to the trace of a solution whose center is an exact frequency-three point.
To apply it one needs even boundary data with nonnegative equatorial trace, close to a normalized odd-degree cone.
The trace must vanish on the equatorial region where the negative conormal of the cone is at least a fixed $\delta>0$.
We obtain this last condition by replacing the even part in a larger ball with the exact Signorini solution and taking its trace on the inner sphere.
Fern\'andez-Real \cite[Theorem 1.1]{FR} constructs an actual thin \emph{free-boundary} point of frequency three in ambient three dimensions; its blow-up is a positive multiple of $|z|^3/3-y_1^2|z|$.
This resolves an existence question for the genuine Signorini equation, not a rate for classical residuals, and its polynomial does not identify the coincidence set of an actual classical solution.

\subsection{Odd error and ordered even comparison}
From \cite[Proposition 3.4, Lemmas 3.6 and 3.8, Propositions 3.9 and 3.10]{FRS} one has $N(r,w)\ge3$, with $N=r\int_{B_r}|Dw|^2/\int_{\partial B_r}w^2$, as well as
\begin{equation}\label{cubic:eq:stock}
 \|w(r\,\cdot)\|_{\mathrm{Lip}(B_3)}\le Cr^3,\quad
 |K\cap B_{2r}|\le Cr^6,\quad
 K\cap B_{2r}\subset\{|z|\le Cr^2\}.
\end{equation}
The constants may depend on the marked pair.
Write $v=w^e$ for the even part.
It is superharmonic, $v|_L\ge0$, and its positive source $\nu=-\Delta v$ is the average of the two reflected volume indicators.
On the support of $\nu$ in $B_{2r}$,
\begin{equation}\label{cubic:eq:trace-support}
 v\le Cr^4 .
\end{equation}
Indeed, on $K$ one side $w=-p_2\le0$, and the other side differs by at most $C r^2|z|\le Cr^4$ by \eqref{cubic:eq:stock} in unscaled coordinates.
The reflected case is identical.

Let $V_{2r}$ be the exact even thin minimizer in $B_{2r}$ with boundary $v$.
The same-boundary comparison and the lower harmonic barrier on each half-ball yield
\begin{equation}\label{cubic:eq:order}
 0\le v-V_{2r}\le Cr^4\quad\text{on }B_{2r}.
\end{equation}
Indeed $V_{2r}\le v$ since $v$ is an admissible superharmonic majorant.
Also $V_{2r}\ge-Cr^2|z|$ by the two half-ball maximum principles, since $v\ge0$ on the plane and $\|Dv\|_\infty\le Cr^2$.
On the support of $\nu$ the difference is then at most $Cr^4$ by \eqref{cubic:eq:trace-support}.
Test its positive part above this bound: the test is zero on the support of $\nu$, and the other solution is superharmonic.
This proves \eqref{cubic:eq:order}.
If $V_r$ instead denotes the same-boundary solution on $B_r$, the same argument and its energy comparison give
\begin{equation}\label{cubic:eq:almost}
 \int_{B_r}|Dv|^2\le \int_{B_r}|D\zeta|^2+Cr^{10}
\end{equation}
for every thin-admissible competitor $\zeta$ of trace $v$.
To see this, first take $\zeta=V_r$, expand its energy about $v$, and bound the difference by $2\int(v-V_r)\,d\nu\le Cr^4 r^6$.
Minimality gives the inequality for every other $\zeta$.

\begin{lemma}[Odd-potential estimate]\label{cubic:le:odd}
On $B_r$,
\begin{equation}\label{cubic:eq:odderror}
 \|w^o-q_o\|_\infty\le Cr^4(1+|\log r|).
\end{equation}
The same remainder bound holds after subtraction of a harmonic Taylor polynomial of degrees at most three for the odd reflected Green potential of $\nu$ on the upper half-space.
\end{lemma}
\begin{proof}
Both potentials are made of pairs of five-dimensional Newtonian kernels at $Y$ and its reflection $Y^*$ with opposite signs, after a smooth cutoff equal to one near zero.
In $B_{2r}$ the support is in an $O(r^2)$ strip and its mass in a ball of radius $s$ is at most $C\min(s^5,r^2s^4)$.
Split sources at distance $Cr^2$ from the evaluation point.
Nearer sources contribute $O(r^4)$ by integrating $|\cdot|^{-3}$ against the bounded densities.
For the farther part of $B_{2r}$, the reflected kernel difference on each shell about the evaluation point is at most $Cr^2s^{-4}$; multiplying by the mass bound gives $Cr^4$ per shell, thus the logarithm.
These estimates can be applied to the pair using the union of the two near balls and then the mean value theorem on the farther part.

On a centered shell of radius $s$, the reflected kernel's Taylor coefficient of order $j$ in the evaluation variable is bounded by $Cs^{-4-j}|z_Y|$.
The mass bound in \eqref{cubic:eq:stock} and $|z_Y|\le Cs^2$ give an absolutely summable series of orders $j\le3$, and the contribution of the near-source Taylor terms on $B_r$ is $O(r^4)$.
For the remaining outer sources the fourth-order Taylor remainder is $O(r^4)$ per centered shell, yielding again $O(r^4(1+|\log r|))$.
All the subtracted Taylor terms are homogeneous harmonic polynomials and odd in $z$.
An ordinary harmonic remainder can be expanded to order three.
The resulting polynomial for $w^o$ must be its known raw cubic $q_o$ by \eqref{cubic:eq:q}; all lower terms vanish.
\end{proof}

\subsection{A power estimate from the odd epiperimetric inequality}
\begin{proposition}\label{cubic:pr:seed}
At each marked point there is $\beta>0$ such that
\begin{equation}\label{cubic:eq:seed}
 \|v-q_e\|_{L^\infty(B_r)}\le C r^{3+\beta}
\end{equation}
for all small $r$.
No neighboring marks or richness assumption is used.
\end{proposition}
\begin{proof}
Use $v_r(X)=r^{-3}v(rX)$ and $c_r=v_r|_{\partial B_1}$.
For the input to \cite[Theorem 1.1]{CV} use instead the trace $\widehat c_r$ on $\partial B_1$ of $\widehat V_r(X)=r^{-3}V_{2r}(rX)$, which is an exact even thin solution in $B_2$.
Both traces are even, nonnegative on the equator, and converge uniformly to $q_e$.
Divide by the fixed positive $\|q_e\|_{L^2(\partial B_1)}$ to put the cone at unit norm in the theorem.
The exact zero trace condition on $\mathcal Z_\delta$ holds for $\widehat c_r$ at every sufficiently small $r$: $\widehat V_r$ converges locally uniformly to $q_e$, and the uniform one-sided $C^{1,1/2}$ interior estimates give one-sided $C^1$ convergence on a neighborhood of the cap.
There the upper normal derivative of the limiting cone is at most a fixed negative constant after the imposed normalization.
Hence the normal derivative of the exact solution is strictly negative there, forcing its thin trace zero by complementarity.
Crucially, we do not impose this condition on $c_r$.
The same nonzero conormal argument is also contained in \cite[Lemma 3.10]{FS}.
No frequency assumption at the center of $\widehat V_r$ is needed by the epiperimetric competitor theorem.

In addition to $0\le c_r-\widehat c_r\le Cr$, we have
\begin{equation}\label{cubic:eq:sphere}
 \|c_r-\widehat c_r\|_{H^1(\partial B_1)}\le Cr^{1/2}.
\end{equation}
Both normalized functions have bounded Lipschitz norms on $B_{3/2}$, by \eqref{cubic:eq:stock} and the exact thin interior estimate.
Their difference is harmonic outside the contact slab $|z|\le Cr$ on the smaller fixed balls.
At a sphere point with $|z|\ge C_1r$ for $C_1$ large enough, a ball of radius $c|z|$ therefore gives gradient bound $Cr/|z|$ by \eqref{cubic:eq:order}.
On the remaining strip use the uniform gradient bound.
Its area is $O(r)$, and the square of the first bound integrates to at most $Cr^2\int_{C_1r}^1 s^{-2}\,ds=O(r)$.
This proves \eqref{cubic:eq:sphere} at each radius without choosing good coarea radii.

Write $Z_c(X)=|X|^3c(X/|X|)$ and $\mathcal W(f)=\int_{B_1}|Df|^2-3\int_{\partial B_1}f^2$.
Let $\widehat\zeta$ be the CV competitor of boundary trace $\widehat c_r$ and contractive factor $1-\kappa\in(0,1)$, after restoring the fixed amplitude.
The value in their proof, with auxiliary exponent $7/2$, is $\kappa=(7/2-3)/(4+7/2+3-1)=1/19$ \cite[(3.7) and the proof of Theorem 1.1]{CV}.
Its $H^1$ seminorm is uniformly bounded by the epiperimetric energy bound, and the same is true of both homogeneous extensions.
Restore the exact actual trace by setting
\begin{equation*}
 \zeta=\widehat\zeta+Z_{c_r-\widehat c_r}.
\end{equation*}
The correction is \emph{nonnegative}, by the exact order in \eqref{cubic:eq:order}, so the competitor is thin-admissible, not just boundary-matched.
Homogeneous extension is bounded from spherical $H^1$ to bulk $H^1$, so \eqref{cubic:eq:sphere} and Cauchy--Schwarz yield
\begin{equation*}
 \mathcal W(\zeta)
 \le(1-\kappa)\mathcal W(Z_{c_r})+Cr^{1/2}.
\end{equation*}
By \eqref{cubic:eq:almost} the rescaled actual energy excess relative to this competitor is $O(r)$, since $r^{5-2+2\cdot3}=r^9$ is the physical Dirichlet energy scale.
The boundary contribution cancels.
Thus, putting $W(r)=\mathcal W(v_r)$,
\begin{equation}\label{cubic:eq:epicompare}
 W(r)\le(1-\kappa)\mathcal W(Z_{c_r})+Cr^{1/2}.
\end{equation}

The following differentiation identity holds, valid for the actual residual without a thin equation,
\begin{equation}\label{cubic:eq:Weiss}
 rW'(r)=9\bigl[\mathcal W(Z_{c_r})-W(r)\bigr]
       +\|\partial_\nu v_r-3c_r\|_{L^2(\partial B_1)}^2.
\end{equation}
For completeness, put $D_r=\int_{B_1}|Dv_r|^2$, $H_r=\int_{\partial B_1}c_r^2$, and $g_r=\partial_\nu v_r$.
Differentiation of the rescaled integrals gives, for almost every $r$,
\begin{align*}
 rD_r'&=-9D_r+\int_{\partial B_1}(|\nabla_\tau c_r|^2+g_r^2),\\
 rH_r'&=2\int_{\partial B_1}c_r(g_r-3c_r).
\end{align*}
Since the radial integral for the degree-three extension is $\int_0^1\rho^8\,d\rho=1/9$, one has $D(Z_{c_r})=(\|\nabla_\tau c_r\|_2^2+9H_r)/9$.
Subtracting $3rH_r'$ from $rD_r'$ and completing the square gives \eqref{cubic:eq:Weiss}.
In particular, this identity uses no thin-obstacle equation for $v$.
We also have a lower bound $W(r)\ge-Cr(1+|\log r|)$.
Indeed the full normalized residual has nonnegative Weiss energy by $N(r,w)\ge3$.
Odd and even pieces are orthogonal.
If $E_o=w^o-q_o$, then $\mathcal W(r^{-3}w^o(r\,\cdot))= \mathcal W(r^{-3}E_o(r\,\cdot))$ exactly: the harmonic degree-three polynomial cancels the cross terms by integration by parts.
Its residual energy is $O(r(1+|\log r|))$ by \eqref{cubic:eq:odderror} and a local squared-cutoff energy test with $E_o$.
In the latter test the total variation of the rescaled signed source is bounded by \eqref{cubic:eq:stock} and the rescaled $E_o$ is $O(r(1+|\log r|))$ in supremum.
The boundary norm is no larger.

Combining \eqref{cubic:eq:epicompare} and \eqref{cubic:eq:Weiss} gives
\begin{equation*}
 r W'(r)\ge c_\kappa W(r)
           +\|\partial_\nu v_r-3c_r\|_2^2-Cr^{1/2},
 \qquad c_\kappa=\frac{9\kappa}{1-\kappa}>0 .
\end{equation*}
Multiplying by $r^{-c_\kappa-1}$ and integrating from $r$ to a fixed small radius gives the upper bound $W(r)\le C r^\vartheta$ for every fixed $0<\vartheta<\min(c_\kappa,1/2)$.
In fact $c_\kappa=1/2$, so one can choose $\vartheta=1/4$.
The lower bound and integration on each $[r/2,r]$ then give
\begin{equation*}
 \int_{r/2}^r
 \|\partial_\nu v_s-3c_s\|_2^2\,\frac{ds}{s}
 \le Cr^\vartheta.
\end{equation*}
Thus by radial differentiation and Cauchy--Schwarz $\|c_r-c_{r'}\|_2\le C r^{\vartheta/2}$ for $r/2\le r'\le r$.
Summing dyadically and using the known limit gives $\|c_r-q_e|_{\partial B_1}\|_2 \le Cr^{\vartheta/2}$ at all small radii.
Integrating such spherical estimates yields the same $L^2(B_2)$ normalized estimate.
The difference of the exact thin solution $\widehat V_r$ and the exact cone $q_e$ has subharmonic square by complementarity \cite[Remark 3.4]{FS}.
The $O(r)$ uniform replacement error and its submean inequality upgrade the bulk estimate to interior supremum, proving \eqref{cubic:eq:seed}, for example with $\beta=1/8$.
\end{proof}

\subsection{Improvement to a logarithmic fourth-order remainder}
Write $A=a e\otimes e$ with $a>0$, let $s=e\cdot y$, and set $\mathcal L(r)=1+|\log r|$ for $0<r<1$.
\begin{proposition}\label{cubic:pr:quartic}
The full classical cubic remainder satisfies
\begin{equation}\label{cubic:eq:quartic}
 \|w-q_e-q_o\|_{L^\infty(B_r)}\le Cr^4\mathcal L(r).
\end{equation}
\end{proposition}
\begin{proof}
It suffices to prove the claim for the even part $v$, since \eqref{cubic:eq:odderror} handles the odd part.
Suppose, starting from Proposition~\ref{cubic:pr:seed}, that $\|v-q_e\|_{L^\infty(B_r)}\le Cr^{3+\beta}$ for some $0<\beta<1$ at all small radii.
On a fixed dilation of $B_r$ compare the ordered exact thin replacement to $q_e$ using \eqref{cubic:eq:order}.
The normalized error is at most $Cr^\beta$.
By the odd case of \cite[Lemma 3.10]{FS} its trace vanishes on the plane away from $|s|\le Cr^{1+\beta/3}$ in $B'_r=B_r\cap L$.
We use its comparison on a larger ball, so that $B'_r$ is strictly interior to the replacement domain.
The rank-one cone is exactly a multiple of the two-dimensional odd model in that lemma.
Thus the actual trace $g=v|_L$ is bounded by $Cr^4$ off that strip, and by $Cr^{3+\beta}$ inside the strip since $q_e|_L=0$.
The strip occupies a fraction at most $Cr^{\beta/3}$ of the four-dimensional disk.
Therefore
\begin{equation}\label{cubic:eq:traceL2}
 \|g\|_{L^2(B'_r)}\le C r^2\bigl(r^4+r^{3+\beta'}\bigr),
 \qquad \beta'=7\beta/6 .
\end{equation}

We next derive an interior estimate from the trace bound \eqref{cubic:eq:traceL2}.
Cut off $g$ and $\nu=-\Delta v$ outside a fixed small neighborhood of the mark.
In the upper half-space use the Poisson extension of this trace and the reflected zero-Dirichlet Green potential of the upper volume source.
Subtracting both from $v$ gives a harmonic function with zero flat trace near the center, whose odd reflection across the plane is analytic.
Lemma~\ref{cubic:le:odd} estimates the volume Green potential after subtracting harmonic Taylor terms through order three by $Cr^4\mathcal L(r)$.
It applies to the upper half of the even source with the same slab and mass bounds.

For the Poisson potential its kernel is $c\,z/(|y-Y|^2+z^2)^{5/2}$, harmonic and odd in $z$ away from the source point.
On each centered tangential annulus of radius $\ell$, its Taylor derivatives of order $j$ at the center are $O(\ell^{-4-j})$.
By \eqref{cubic:eq:traceL2} the $L^1$ trace mass there is at most $C(\ell^8+\ell^{7+\beta'})$.
Thus its harmonic Taylor coefficients through order three are integrable.
For sources $|Y|\le2r$ the Poisson extension on $B_r^+$ has $L^2$ norm at most $Cr^{1/2}\|g\|_{L^2(B'_{2r})}$ by the $L^1$ normalization of the Poisson kernel on each plane at height $z$ and integration in $z$.
Dividing by the square root of the bulk volume gives $C(r^4+r^{3+\beta'})$.
The Taylor terms of those sources have at most the same normalized order, by summing $Cr^j(\ell^{4-j}+\ell^{3+\beta'-j})$ over $\ell\lesssim r$ for $j\le3$.
For the remaining Poisson sources, the fourth-order Taylor remainder on a shell of radius $\ell\ge2r$ is bounded by $Cr^4(1+\ell^{-1+\beta'})$.
Summing shells yields
\begin{equation*}
 C\bigl(r^4\mathcal L(r)+r^{3+\beta'}\bigr)
\end{equation*}
as a normalized bulk $L^2$ error, including the case $\beta'=1$.
Hence $v$ in the upper half-space is a harmonic polynomial of degree at most three plus this error.
The lower-degree coefficients are zero and the cubic equals $q_e$ restricted to the upper half-space by the raw limit.
Even reflection gives the same estimate in the full ball.
Again use the submean inequality for the square of the difference of the exact thin replacement and $q_e$ on a larger ball, together with the uniform $O(r^4)$ replacement error.
The normalized bulk estimate upgrades to
\begin{equation*}
 \|v-q_e\|_{L^\infty(B_r)}
 \le C\bigl(r^4\mathcal L(r)+r^{3+\beta'}\bigr).
\end{equation*}
While $\beta'<1$ this improves $\beta$ to $7\beta/6$.
After finitely many iterations one has $\beta'\ge1$, and the last display gives the claimed single logarithm.
No free-boundary regularity beyond the thin strip lemma has been assumed.
\end{proof}

\subsection{The backward parameter estimate}
\begin{theorem}\label{cubic:th:param}
The projection onto the parameter axis of the five-dimensional rank-one nonharmonic cubic cohort is Lebesgue null under the monotonicity and uniform boundary increase of \cite{FRS}.
On localized pieces it has the pairwise bound
\begin{equation}\label{cubic:eq:time}
 |t_1-t_0|\le C|x_1-x_0|^3\mathcal L(|x_1-x_0|).
\end{equation}
\end{theorem}
\begin{proof}
Choose the later mark as base $(x_0,t_0)=(0,0)$; the other mark has $t_1\le0$.
Set $r=|x_1|>0$, and use the base normal $z=X\cdot\nu_0$.
Work in a compact interior cylinder.
The normals of nearby quadratic marks can be aligned with $|\nu_1-\nu_0|\le Cr$ and their center offsets obey $|x_1\cdot\nu_i|\le Cr^2$ for $i=0,1$.
These rates follow directly from the bounded normalized quadratic comparison polynomial in the proof of \cite[Lemma 6.5, (6.25)]{FRS}.
In detail, the polynomial
\begin{equation*}
 r^{-3}\bigl[p_{2,0}(x_1+rX)-p_{2,1}(rX)\bigr]
\end{equation*}
is bounded on a fixed ball for nearby same-cohort pairs at the fixed later base: in the three-way normalization of that proof the scale denominator is comparable to $r^3$, since the translated old nonzero raw cubic survives with positive coefficient.
Comparing the degree-two coefficients gives $|\nu_1\otimes\nu_1-\nu_0\otimes\nu_0|\le Cr$, and the degree-one coefficient gives $|x_1\cdot\nu_0|\le Cr^2$.
The other normal offset follows from the first angle bound.
A contradiction sequence proves both bounds for all sufficiently close pairs at a fixed later base.

Define the harmonic jet subtractions $P_i=p_{2,i}+q_{o,i}$ in the corresponding centered ambient coordinates.
At $x_1+rX$, the earlier residual is equal to the translated later residual plus
\begin{equation*}
 T(X)=u^{t_1}(x_1+rX)-u^0(x_1+rX)\le0
\end{equation*}
and the polynomial $J(X)=P_0(x_1+rX)-P_1(rX)$.
The polynomial $J$ is harmonic, since both second Laplacians are the constant $r^2$ in $X$.
On the later plane $L_0$ its trace is $O(r^4)$ on fixed balls by the normal angle and shift bounds: the quadratic uses a normal coordinate $O(r^2)$, and the odd cubic is that coordinate times a quadratic of size $O(r^2)$.
Hence its normally even part is $O(r^4)$, by finite-dimensional norm equivalence and injectivity of flat trace on even harmonic polynomials.

Use \eqref{cubic:eq:quartic} for the two actual residuals after $P_i$ is subtracted.
Rotating the new even cubic to the old plane, and replacing $x_1/r$ by its projection $Y$ on the old plane in the old translated cubic, cost only $O(r^4)$ in physical values by homogeneity, Lipschitz continuity on fixed balls, and the just established angle and shift bounds.
On the common unit sphere the resulting even equation therefore has the form
\begin{equation*}
 r^3 q_{e,1}(X)=
 r^3 q_{e,0}(X+Y)+T^e(X)+O(r^4\mathcal L(r)),
\end{equation*}
with both $q_{e,i}$ understood on the common plane.
Take its inner product with the positive function $|z|$.
Each centered cubic is orthogonal to it: on the upper half-sphere the restrictions are Dirichlet spherical eigenfunctions of degrees three and one, respectively, and even reflection doubles the inner product.
The translated cubic satisfies
\begin{equation*}
 q_{e,0}(X+Y)-q_{e,0}(X)
  =-|z|\bigl(2y^{\mathsf T}A_0Y+Y^{\mathsf T}A_0Y\bigr).
\end{equation*}
The first term integrates to zero since $y\mapsto-y$ preserves both the sphere and the weight $z^2$.
If $M_2=\int_{\partial B_1}z^2$, the weighted identity therefore reads
\begin{equation*}
 0=-r^3 M_2\,Y^{\mathsf T}A_0Y
     +\int_{\partial B_1}|z|T^e+O(r^4\mathcal L(r)).
\end{equation*}
The two displayed terms on the right are nonpositive: $A_0$ is positive semidefinite and $T^e\le0$ by the order of the leaves.
Hence $\int_{\partial B_1}|z|\,|T^e|\le Cr^4\mathcal L(r)$.
Reflection preserves $|z|$, and $T$ itself is nonpositive, so the weighted integrals of $|T^e|$ and $|T|$ are equal.
This ordering of the two actual leaves is essential: the proof does not estimate the sign-opposite forward pair in the same coordinates.

By \cite[Lemma 9.2(b)]{FRS} applied at the earlier actual mark, the time increase $-T$ is at least $c\,r|t_1|$ on its two transverse caps on this same comparison sphere.
Aligning the two normals by $O(r)$ ensures $|z|$ is bounded below on a fixed subcap of positive area.
This proves \eqref{cubic:eq:time}.
All residual constants, radii and earlier warm comparison constants can be made uniform by countably many localizations.
Past exclusion \cite[Lemma 9.4]{FRS} prevents two distinct marked times over the same spatial point (if needed, use same-point stationarity of \cite[Lemma 6.2]{FRS} on the interval).
Compactness of each interior parameter interval and closedness of the full singular set give local approach in time as well as in space.
Thus the pairwise estimates apply on sufficiently small spatial pieces.

Finally use the countably many three-dimensional spatial Lipschitz charts from \cite[Proposition 8.3]{FRS}, based on the secant estimate in \cite[Lemma 6.6]{FRS}.
Pull back \eqref{cubic:eq:time} to each chart domain; its inverse and forward Lipschitz constants only modify $C$ and the small radius.
The critical logarithmic Sard bound for subsets of $\R^3$ with modulus $Cr^3\mathcal L(r)^2$ (proved in Theorem~\ref{sard:th:sard}) applies with room to spare.
This proves that the parameter projection is null on each piece and hence on their countable union.
\end{proof}

\section{The harmonic cubic with next frequency \texorpdfstring{$7/2$}{7/2}}\label{sec:slit}

\subsection{Rescalings at neighboring singular points}
We work on a countably localized piece $E$ of the quadratic-rank-one classical obstacle stratum, with completed harmonic cubic and next frequency $\lambda=7/2$, for the monotone families of \cite{FRS}.
Write $w_x=u^{t_x}(x+\cdot)-P_x$ for $P_x=p_{2,x}+p_{3,x}$, $L_x=\ker p_{2,x}$, and $h_x(r)=\|w_x(r\,\cdot)\|_{L^2(\partial B_1)}$.
One can partition by integer bounds for the interior cylinders, the size of the cubic jets, the constants and the small-radius thresholds in the local estimates used below.
In particular fix $\eta=1/16$ and use
\begin{equation}\label{slit:eq:stocks}
 cr^{\lambda+\eta}\le h_x(r)\le Cr^\lambda,\quad
 \|w_x(r\,\cdot)\|_{\mathrm{Lip}(B_3)}\le Ch_x(r),\quad
 N(r,w_x)\ge\lambda-Cr^{1/4}.
\end{equation}
Here $N(r,f)=r\int_{B_r}|Df|^2/\int_{\partial B_r}f^2$.
These follow from \cite[Lemmas 4.1, 4.3, 6.7 and Proposition 4.12(a)]{FRS}.
For clarity, the one-sided frequency deficit follows directly from \cite[Lemma 4.3]{FRS} with $\gamma=15/4$.
Writing $H_x=h_x(r)^2$ and $D_x=r^2\int_{B_1}|Dw_x(rX)|^2\,dX$, define the truncated frequency
\begin{equation*}
 \phi_\gamma(r,w_x)=\frac{D_x+\gamma r^{2\gamma}}{H_x+r^{2\gamma}}.
\end{equation*}
The cited differential inequality gives $\phi_\gamma'(r,w_x)\ge-Cr^{3-\gamma}$; its limit is $\lambda$, so integration yields $\phi_\gamma(r,w_x)\ge\lambda-Cr^{1/4}$.
Moreover,
\begin{equation*}
 N(r,w_x)=\phi_\gamma(r,w_x)
 +(\phi_\gamma(r,w_x)-\gamma)\frac{r^{2\gamma}}{H_x},
 \qquad \frac{r^{2\gamma}}{H_x}\le Cr^{3/8}.
\end{equation*}
Since $\phi_\gamma$ is locally bounded, this proves the last inequality in \eqref{slit:eq:stocks}.
The ordinary and completed cubic subtraction differ by $O(r^4)$.
On the contact, $\nabla u=0$ and $P_x=z^2/2+za_x$, where $a_x$ is quadratic in the normal coordinate $z$ and the tangential coordinates.
The gradient bound in \eqref{slit:eq:stocks} gives $|z|\le Cr^2$ and $|z+a_x|\le Cr^{5/2}$ on contact in $B_r$.

For each fixed dilation factor $\Lambda\ge1$, two power estimates are important here:
\begin{equation}\label{slit:eq:oddreplacement}
 \|w_x^o\|_{L^\infty(B_r)}\le Cr^4,\qquad
 \|v-w_x^e\|_{L^\infty(B_{2\Lambda r})}\le C_\Lambda r^4 .
\end{equation}
The superscripts denote odd and even normal reflection, and $v$ is the even exact thin minimizer on a ball of radius $2\Lambda r$ centered at this mark with trace $w_x^e$ and plane $L_x$.
For the first assertion we use reflected Newtonian potentials.
Indeed in $B_{2r}$ the contact has thickness $O(r^{5/2})$ on normal lines and normal position $O(r^2)$.
The reflected pair of Newtonian sources for the odd part gives $O(r^{9/2})$ from distances $\le Cr^2$ (integrate the mass bound $C\min(s^5,r^{5/2}s^4)$ against the Newtonian kernel); at greater distances the paired kernel is $O(r^2s^{-4})$ and the same mass bound gives $O(r^{9/2})$ per shell.
Thus the near-source bound is $O(r^{9/2}(1+|\log r|))$.
The coefficients through order four at the center for the paired potential are absolutely integrable: on centered shells at $s$, the paired derivatives of order $j$ are $O(s^{-4-j}s^2)$ and the mass is $O(s^{13/2})$.
The latter also follows by cutoff from $-\Delta w_x={\bf1}_{K_x}$ and \eqref{slit:eq:stocks}.
A fifth-order Taylor remainder for the farther sources is $O(r^{9/2})$.
The analytic harmonic part and these Taylor terms cancel in degrees below four because $w_x=O(r^{7/2})$.
Constants are uniform on localized pieces.

For the second assertion first use the exact thin solution with the full boundary datum $w_x$.
Its upper bound is $w_x$ by superharmonic comparison, and its lower bound on either half is $-C_\Lambda r^{5/2}|z|$ by the harmonic comparison with this function.
On actual contact $w_x=-P_x\le C_\Lambda r^4$ by completing the square.
Test the positive part of the difference above $C_\Lambda r^4$: it vanishes on the volume source of $w_x$ and the thin solution is superharmonic.
This gives the uniform error; its even part is again exact and has trace $w_x^e$.
We use this construction on each fixed dilated ball, not a single global thin approximant.

Fix the old base $x=0$, $t_x=0$, normal $\nu_0$, and consider matched actual neighbors $x_j,t_j$ with $|x_j|\le r/2$ in $E$ along $r\downarrow0$.
Set $h=h_0(r)$ and
\begin{align*}
 Y_j&=x_j/r,\qquad F_0(X)=w_0(rX)/h,\quad F_j(X)=w_{x_j}(rX)/h,\\
 T_j(X)&=(u^{t_j}(x_j+rX)-u^0(x_j+rX))/h,\\
 J_j(X)&=(P_0(x_j+rX)-P_{x_j}(rX))/h .
\end{align*}
Thus $F_j=F_0(Y_j+\cdot)+T_j+J_j$.
On a subsequence with $F_0\to q$ and $Y_j\to Y_j^*$, \cite[Lemma 6.8 and its proof]{FRS} give
\begin{equation}\label{slit:eq:match}
 q(Y_j^*+\cdot)=q,\quad F_j\to q,\quad T_j,J_j\to0
 \quad\text{locally uniformly}.
\end{equation}
The required strong $W^{1,2}_{\rm loc}$ convergence of the base follows from the normalized Lipschitz bounds and harmonic interior gradient convergence away from the vanishing contact slab: omit a fixed thin strip first and then let its width tend to zero.
For the termwise consequences, divide by the sum of the norms of the three terms as in that proof.
After the invariance is established the time and jet limits sum to zero.
The time limit has one sign, and the harmonic jet limit has a zero: the two ends of the centered segment $[-Y_j,0]$ give opposite weak signs up to $O(r^4)$ before normalization.
These limits therefore vanish.
The translated base term is asymptotically of norm $h$, so the summed normalization is comparable to $h$.

Comparing the quadratic coefficients of $J_j$ gives
\begin{equation*}
 |\nu_j\otimes\nu_j-\nu_0\otimes\nu_0|
 \le Cr+C(h/r^2)\|J_j\|_{L^2(\partial B_1)}\le Cr.
\end{equation*}
The translated cubic contributed $O(r^3)$ before division.
Take the sign of $\nu_j$ close to $\nu_0$.
By time order either $x_j\in K_0$ or $0\in K_j$; the contact slab then gives $|x_j\cdot\nu_0|+|x_j\cdot\nu_j|\le Cr^2$ in both cases.
Consequently $J_j|_{L_0}=O_\Lambda(r^4/h)$ on fixed balls.
Finite dimensionality and injectivity of restriction to $L_0$ on the even harmonic polynomials of degree at most three show that $J_j^e=O_\Lambda(r^{7/16})$.
Rotating the neighbor's plane to the base plane costs $O_\Lambda(r)$ using the normalized Lipschitz bounds, as does replacing $Y_j$ by its tangential projection $\bar Y_j$.
Denote the actual exact even thin replacements, after this alignment, by $V_0,V_j$.
We have
\begin{equation}\label{slit:eq:matchedmodel}
 V_j=V_0(\bar Y_j+\cdot)+T_j^e+O_\Lambda(r^{7/16})
\end{equation}
on the controlled interior of each fixed dilated ball.

For any fixed $b>0$ the models have
\begin{equation}\label{slit:eq:freq}
 N(s,V_i)\ge\lambda-C_b r^{1/8}\quad(b\le s\le1/2).
\end{equation}
Indeed the full original residual and its even thin replacement differ in supremum norm by $O_\Lambda(r^{7/16})$ by \eqref{slit:eq:oddreplacement}.
The two nonnegative source measures have bounded mass locally by cutoff and the supremum bounds.
Testing their difference with its product with a squared cutoff gives gradient $L^2$ error $O_\Lambda(r^{7/32})$.
The height denominators on fixed positive scales stay bounded below by \eqref{slit:eq:match}, so \eqref{slit:eq:stocks} gives \eqref{slit:eq:freq}.
The model's exact frequency at radius zero is not used.

In the compact interior parameter intervals used here, spatial convergence of marked points in the frequency-$7/2$ set to a fixed spatial coordinate of the same set also yields convergence of the parameters.
Indeed, a subsequential limit is a singular pair by \cite[Lemma 6.2(a)]{FRS}.
If its parameter were different, \cite[Lemma 6.2(c)]{FRS} would give local stationarity between the two parameters.
At the earlier pair the jets and frequency would therefore be the same as those at the prescribed pair.
The forward estimate following from \eqref{eq:comparison} and \cite[Lemma 9.2(b)]{FRS} excludes such stationarity: at the very same spatial coordinate the error is $O(r^{7/2})$, while the lower bound for a fixed positive parameter increment is $cr(t-t_0)$.
Letting $r$ tend to zero gives strict positivity there.
The same argument justifies the use of spatial secants in the frequency-$7/2$ chart construction below.

\subsection{The limit of the normalized parameter increment}
\begin{lemma}\label{slit:le:time}
Assume the matched limit $q$ in \eqref{slit:eq:match} is a pure slit.
Use coordinates $(y,s,z)$ for its three passive and two slit variables, and write $U_\kappa(s,z)= \operatorname{Re}(s+i|z|)^\kappa$ with argument in $[0,\pi]$.
If $m_k\ge r_k^\delta$ for $0<\delta<7/16$ and $m_k\ge\|T_{j,k}^e\|_{L^\infty(B_{3/4})}$, then, along a subsequence, $T_{j,k}^e/m_k$ converges in $L^1_{\rm loc}$ to $t_j U_{1/2}$ for a real $t_j$.
If the second inequality is an equality, the limit is nonzero.
\end{lemma}
\begin{proof}
Pass to one sign of the time increment.
The absolute increment $f_k=|T_{j,k}^e|/m_k$ is subharmonic: the ordered difference of solutions and their positivity indicators have the same weak sign.
Its Laplacian is supported in the two contact slabs and their reflections, within $|z|\le C_\Lambda r_k$ on each fixed ball $B_\Lambda$.
This follows from the normal alignment and translations above.
By \eqref{slit:eq:match} and the lower bound on $m_k$, we have $f_k\le C_\Lambda r_k^{-\delta}$ on each fixed ball.
Hence the $L^1$ integral in this shrinking slab is $O_\Lambda(r_k^{1-\delta})$.

Outside the slab $f_k$ is nonnegative harmonic.
Its values at fixed warm heights on either side of the slab are bounded by Harnack from $B_{3/4}$, where its supremum is at most one.
On each horizontal section above the slab, the integral over a fixed tangential disk is bounded by a constant times one such warm value.
To see this, take a half-ball directly above the section with its flat disk substantially larger than the target disk; discard the nonnegative curved boundary data and use the positive Poisson kernel at an interior point a fixed height above the section, bounded below on the target disk.
The constant is uniform in the height of the sections in a fixed compact set down to the edge of the slab.
Apply the same argument below.
This proves uniform local $L^1$; the subharmonic submean inequality then gives uniform local supremum bounds over the whole ambient space.
Extract a locally harmonic limit off $L_0$; local dominated convergence gives $L^1_{\rm loc}$ convergence.
If the supremum over $B_{3/4}$ is one, the submean inequality on a slightly larger fixed ball gives a positive lower bound on the $L^1$ norm there, hence the limit is nonzero.

For any fixed compact patch of $\{z=0,s>0\}$, both exact states in \eqref{slit:eq:matchedmodel} have positive thin trace for all sufficiently large $k$ by \eqref{slit:eq:match}.
They are harmonic across that patch.
The error $O_\Lambda(r_k^{7/16})/m_k$ tends to zero, so the limit is harmonic there and has zero upper normal derivative by evenness.
On a compact patch of $\{z=0,s<0\}$, both exact states have zero trace for all large $k$ by \cite[Lemma 3.9]{FS}, or by one-sided $C^1$ compactness and the strictly negative upper derivative of the pure limit.
Their rescaled difference is bounded by \eqref{slit:eq:matchedmodel} and the already proved local supremum bounds.
Its upper harmonic restriction extends oddly across the flat patch, so that its limit has zero trace.
These arguments apply on every fixed compact patch; they do not assume one thin approximation on the whole space.

The nonnegative limit is a harmonic function with mixed boundary conditions on the upper half-space, with Dirichlet zero on the negative slit and Neumann zero on the opening.
The codimension-two edge does not generate a boundary term: write $\rho=(s^2+z^2)^{1/2}$ and let $\chi_\varepsilon$ be zero for $\rho\le\varepsilon^2$, one for $\rho\ge\varepsilon$, and affine in $\log\rho$ between the two radii.
For each bounded set the integral of $|D\chi_\varepsilon|^2$ is $O(|\log\varepsilon|^{-1})$, since the transverse area element is $\rho\,d\rho\,d\theta$.
With a compact cutoff $\eta$, test the weak mixed equation away from the edge with $\eta^2\chi_\varepsilon^2 f$, where $f$ denotes the limit.
The Dirichlet trace of the test is zero.
Cauchy--Schwarz bounds $\int\eta^2\chi_\varepsilon^2|Df|^2$ by $C\|f\|_\infty^2\int(|D\eta|^2+\eta^2|D\chi_\varepsilon|^2)$.
Letting $\varepsilon\downarrow0$ gives finite energy up to the edge.
For a smooth mixed test $\varphi$, one may now use $\chi_\varepsilon\varphi$ and pass to the limit.
The extra term $\int\varphi Df\cdot D\chi_\varepsilon$ tends to zero by Cauchy--Schwarz and the same cutoff bound.
The separated mixed spherical degrees are $\kappa=m+1/2+\ell+2j$, with $m,\ell,j\ge0$ integral, as in the half-slit spectral decomposition in \cite[Section 3]{FS}.
In transverse polar coordinates the factors in the transverse angle are $\cos((m+1/2)\theta)$; the remaining spherical factors are bounded passive harmonics and Jacobi polynomials multiplied by nonnegative powers of the two radial lengths.
In particular the absolute value of each eigenfunction is bounded by a constant times the ground eigenfunction $U_{1/2}|_{\partial B_1}$: the quotient of the two cosine factors is bounded up to $\theta=\pi$, and the transverse radial power is at least $1/2$.
For each spherical mode its radial coefficient is a multiple of $R^\kappa$ by the harmonic radial ODE; local boundedness excludes the singular solution.
Positivity bounds the absolute value of every coefficient by a constant times the ground coefficient $cR^{1/2}$.
Letting $R\to\infty$ eliminates all higher modes.
Restore the sign of the increment.
\end{proof}

\subsection{A height estimate with an additive error}
The sliding maximum argument below is adapted from \cite[Proposition 5.2]{FS}.
Their calculation is for an exact Signorini solution with prescribed frequency at the center.
Here the comparison solutions have an error in the frequency bound on each interval of positive radii.
We retain that error throughout the calculation.
Let $V_k$ be one of the exact thin models, and $q_k$ a pure slit of positive amplitude bounded away from zero, with $Z_k=V_k-q_k$ tending to zero locally uniformly.
With $\mu=\lambda-1/3$, write $H(s,f)=\|f(s\,\cdot)\|_{L^2(\partial B_1)}^2$.
For each fixed $a>0$, \eqref{slit:eq:freq} and two Green identities give on $[a,1/2]$
\begin{equation}\label{slit:eq:deriv}
 \frac{d}{ds}[s^{-2\mu}H(s,Z_k)]
 \ge \frac{2s^{-2\mu}}s\left[
 (\lambda-\mu)H(s,Z_k)-\int_{B_1}(Z_k)_s\Delta(Z_k)_s
 -C_a r_k^{1/8}H(s,V_k)\right].
\end{equation}
In fact $Z_k\Delta Z_k=V_k(-\Delta q_k)+q_k(-\Delta V_k)\ge0$.
Expanding the gradient square and using homogeneity of $q_k$ for the boundary cross term gives $D((Z_k)_s)\ge\lambda H(s,Z_k) -2\int(Z_k)_s\Delta(Z_k)_s-C_a r_k^{1/8}H(s,V_k)$, and the height derivative gives \eqref{slit:eq:deriv}.
Take the sliding maximum of $s^{-2\mu}H(s,Z_k)$ over $[s,2s]$.
If $s$ itself maximizes $G_k(\sigma)=\sigma^{-2\mu}H(\sigma,Z_k)$ on $[s,2s]$, then for $1\le \rho\le2$, $H(\rho s,Z_k)\le\rho^{2\mu}H(s,Z_k)$.
Hence
\begin{equation*}
 \|(Z_k)_s\|_{L^2(B_2\setminus B_1)}^2
   =\int_1^2\rho^4H(\rho s,Z_k)\,d\rho
   \le H(s,Z_k)\int_1^2\rho^{4+2\mu}\,d\rho .
\end{equation*}
The superquadratic bound \cite[Lemma 3.12]{FS} then bounds the middle term by $o(1)H(s,Z_k)$ uniformly on fixed $[a,b]$, $b<1/4$, since the normalized slit amplitude on that interval stays positive.
If a maximizer lies strictly to the right, it persists in the sliding interval for small positive steps, so the maximum has nonnegative right lower derivative.
The maximum is locally Lipschitz; integrating the resulting lower derivative bound gives
\begin{equation}\label{slit:eq:height}
 s^{-2\mu}H(s,Z_k)\le
 C_b\|Z_k\|_{L^\infty(B_{2b})}^2+C_a r_k^{1/8}
 \quad(a\le s\le b).
\end{equation}
The first constant does not depend on $a$.
For any denominator $m_k\ge r_k^{1/32}$, the additive error vanishes after division by $m_k^2$.
Only fixed positive scales $s$ are sent through the thin model limit; its exact limiting frequency at radius zero is irrelevant.

\subsection{One-step improvement at three-direction scales}
At the fixed old base define
\begin{equation*}
 {\cal A}(r)=\inf_{a\ge0,\ e\in L_0,\ |e|=1}
 \|w_0^e(r\,\cdot)-a U_\lambda(e\cdot X,z)\|_{L^2(B_1)},
 \qquad d(r)={\cal A}(r)/h_0(r).
\end{equation*}
We use the norm on the ball because the identities obtained by differentiating at a minimizing profile pass to the limit under weak convergence in $L^2(B_1)$.
Call a scale $\epsilon$-rich when the normalized positions of same-cohort neighbors of $E$ in $B_{r/2}(0)$ cannot all be covered by an $\epsilon$-tube about any linear two-plane.
At such a scale one can choose three of them with successive distances to the spans of the preceding ones bounded below in terms of $\epsilon$, by the greedy selection of independent vectors.
This is a uniform linear independence after taking limits.

\begin{theorem}\label{slit:th:one}
For the fixed base of $E$ and any $\epsilon>0$ there is $\sigma>0$ such that for each $\gamma>0$ one can choose a threshold $\eta_0>0$ so that for all sufficiently small $r$ with $d(r)\le\eta_0$,
\begin{equation}\label{slit:eq:dyad}
 \frac{{\cal A}(r/2)}{h_0(r)}
 \le
 \begin{cases}
 2^{-\lambda-\sigma}d(r)+C r^{1/32},
       &r\text{ is }\epsilon\text{-rich},\\
 2^{-\lambda+\gamma}d(r)+C r^{1/32},
       &\text{otherwise}.
 \end{cases}
\end{equation}
The threshold and constants may depend on the base, localized bounds, and the specified parameters.
On all sufficiently small rich scales $d(r)$ is small.
\end{theorem}
\begin{proof}
At a sequence of rich scales any normalized limit has the three independent limiting directions in its translation spine by \cite[Lemma 6.8]{FRS}.
The two-dimensional classification \cite[Lemma A.2]{FRS} forces a pure slit, so $d(r)\to0$ along such scales.

To disprove the near-nonworsening inequality one could choose $r_k\downarrow0$, $d_k=d(r_k)\to0$ and $d_k/r_k^{1/32}\to\infty$, with ratio of the inner and outer errors greater than $2^{-\lambda+\gamma}$.
Indeed if not, there are thresholds and an additive constant as stated; if $d_k/r_k^{1/32}$ were bounded, the inner error by restriction of the outer optimizer would satisfy the floor bound.
For strict gain, if no positive strict factor exists, take likewise gains tending to zero, with outgoing ratio larger than $2^{-\lambda-\sigma_k}$ for $\sigma_k\downarrow0$.
Choose optimal $q_k$ after dividing the outer comparison by $h_0(r_k)$.
An exact even thin replacement $V_{0,k}$ on a ball containing $\overline B_1$ gives $Z_{0,k}=V_{0,k}-q_k$ with $\|Z_{0,k}\|_{L^2(B_1)}/d_k\to1$ by \eqref{slit:eq:oddreplacement}.
The square $Z_{0,k}^2$ is subharmonic.
Submean and Caccioppoli give local uniform supremum and energy bounds on $Z_{0,k}/d_k$.
Take its weak $L^2(B_1)$ and strong interior $L^2$ limit $z_*$.
Both phases are fixed off the vanishing edge strip by \cite[Lemma 3.9]{FS}, so the limit is a fixed-slit mixed solution, with the weak edge justification as in Lemma~\ref{slit:le:time}.
Estimate \eqref{slit:eq:height} gives $H(s,z_*)^{1/2}\le Cs^\mu$.
The mixed degrees are half-integers, so all modes below $\lambda$ are absent.
By orthogonality on balls,
\begin{equation}\label{slit:eq:limitnorm}
 \|z_*(\cdot/2)\|_{L^2(B_1)}
 \le2^{-\lambda}\|z_*\|_{L^2(B_1)}
 \le2^{-\lambda}.
\end{equation}
The inner optimal error divided by the outer error is at most the corresponding limit on the left.
This proves the near-nonworsening inequality in contradiction form.

For the strict inequality, choose three actual matched neighbors on each rich scale.
Their limiting positions $Y_j^*$ span the passive edge of the limiting slit.
Let $b_{j,k}=e_k\cdot\bar Y_{j,k}$ be the opening offset of neighbor $j$ for the optimal orientation $e_k$; it tends to zero.
Use the denominator
\begin{equation*}
 m_{j,k}=\max\{d_k,|b_{j,k}|,\|T_{j,k}^e\|_{L^\infty(B_{3/4})}\}.
\end{equation*}
For fixed small $b>0$ the centered difference $Z_{j,k}=V_{j,k}-q_k$ is bounded on $B_{2b}$ by $Cm_{j,k}$.
Indeed use \eqref{slit:eq:matchedmodel}, the outer error bound on a compact subball of $B_1$ containing $\bar Y_{j,k}+B_{2b}$, the change of the homogeneous slit due to the offset, and the time normalization.
Choose $b<1/8$, for instance.
If $m_{j,k}/d_k\to\infty$, the centered difference divided by $m_{j,k}$ has limit $c U_{5/2}+tU_{1/2}$: the base error vanishes on the translated interior, the slit offset differentiates the pure profile in its opening coordinate, and Lemma~\ref{slit:le:time} gives the time mode.
At least one coefficient is nonzero by the choice of denominator and the nonvanishing assertion for the time mode.
Estimate \eqref{slit:eq:height} gives growth $O(s^\mu)$, impossible for these two distinct degrees.
It follows that $m_{j,k}=O(d_k)$.
Reapplying the estimate at denominator $d_k$ gives the limit bound
\begin{equation}\label{slit:eq:translation}
 H(s,z_*(Y_j^*+\cdot)+c_jU_{5/2}+t_jU_{1/2})^{1/2}
 \le C s^\mu .
\end{equation}
Sphere convergence at each fixed $s$ follows by local uniform convergence off the thin plane and uniform supremum bounds, since the plane cuts a surface-null subset of the sphere.
For the base normalized error one uses interior mixed harmonic compactness as well.
Both added modes are independent of the passive variables.
Caccioppoli for the mixed solution in \eqref{slit:eq:translation}, integrating the spherical bound over the smaller ball, yields
\begin{equation*}
 s\|\partial_{y_i}z_*(Y_j^*+s\,\cdot)\|_{L^2(B_{1/2})}
 \le C s^\mu,\qquad i,j=1,2,3 .
\end{equation*}

In a strict-gain contradiction the limiting inner error is at least $2^{-\lambda}$.
Both inequalities in \eqref{slit:eq:limitnorm} are therefore equalities: $\|z_*\|_{L^2(B_1)}=1$, and every mixed mode of degree strictly larger than $\lambda$ vanishes.
Thus $z_*$ is a nonzero $\lambda$-homogeneous mixed solution.
Writing $\rho^2=s^2+z^2$, its separated degree-$7/2$ form is
\begin{equation*}
 c U_{7/2}+L_1(y)U_{5/2}
 +[Q_2(y)-(\Delta_y Q_2/10)\rho^2]U_{3/2}
 +[P_3(y)-(\Delta_y P_3/6)\rho^2]U_{1/2},
\end{equation*}
where the passive polynomials have degrees $1,2,3$.
The coefficients follow from $\Delta_{s,z}(\rho^2U_\alpha)=4(\alpha+1)U_\alpha$.
This spans the separated degrees of Lemma~\ref{slit:le:time}.
After differentiation in $y_i$ and translation by $Y_j^*$, the degree-$3/2$ terms are
\begin{equation*}
 \partial_iQ_2(Y_j^*)U_{3/2}
 \quad\hbox{and}\quad
 (D^2P_3(Y_j^*)y)_iU_{1/2}.
\end{equation*}
Their different transverse angular factors make them independent.
Since $\mu-1=13/6>3/2$, the derivative growth bound forces both coefficients to vanish; terms of degree $1/2$ must vanish as well.
The maps $\nabla Q_2$ and $D^2P_3$ are linear and vanish at the three independent vectors $Y_j^*$, so both vanish identically.
Homogeneity then gives $Q_2=P_3=0$.
The remaining amplitude and rotation modes are orthogonal to $z_*$ by differentiating the original ball-norm optimality of $q_k$ and sending $k\to\infty$ by weak $L^2(B_1)$.
This is a contradiction.
\end{proof}

\subsection{Approximation of the coincidence set by a slit}\label{sec:slit-contact}

Write $x=(y,z)\in\R^4\times\R$, let $r=|x|$, and let $K=\{u^0=0\}$ at the base time.
Put $w=u^0-\mathcal P$, with $\mathcal P$ the completed harmonic cubic ansatz of \cite[Definition 4.5]{FRS}; in particular $\Delta\mathcal P=1$.
At frequency $7/2$, \cite[Lemma 4.1(b), Lemma 4.7, Lemma 4.9 and Proposition 4.12(a)]{FRS} give
\begin{equation}\label{access:eq:amp}
 \|w\|_{L^\infty(B_r)}+r\|Dw\|_{L^\infty(B_r)}\le Cr^{7/2},
 \qquad H(r,w)^{1/2}\ge c_\epsilon r^{7/2+\epsilon}
\end{equation}
for sufficiently small $r$ and each fixed small $\epsilon>0$.
Indeed the lower frequency ratio in \cite[Lemma 4.1(b)]{FRS} gives the upper $H(r,w)\le Cr^7$ using any admissible truncation order above four; the reverse bound with an arbitrarily small power loss follows from \cite[Lemma 4.1(a)]{FRS}.
\cite[Lemma 4.7]{FRS} provides the gradients by the annular $L^2$ bound and $w(0)=0$ gives the function bound.

On $K$, $Du^0=0$.
Write $a(y,z)=p_3(y,z)/z$ as a quadratic polynomial.
The ansatz is $\mathcal P=(z+a)^2/2+zQ$, with $Q$ a fixed homogeneous cubic.
Hence
\begin{equation*}
 \partial_z\mathcal P=(z+a)(1+\partial_z a)+O(r^3).
\end{equation*}
Using \eqref{access:eq:amp} and absorbing the bounded perturbation of $1$ in the factor, one obtains
\begin{equation}\label{access:eq:vert}
 |z+a(y,z)|\le Cr^{5/2},\qquad |z|\le Cr^2
 \quad\hbox{on }K\cap B_r .
\end{equation}
On this set $\mathcal P=O(r^5)$, since the square is $O(r^5)$ and $zQ=O(r^5)$.
Accordingly
\begin{equation}\label{access:eq:tiptrace}
  \sup_{K\cap B_r}|w|/H(r,w)^{1/2}\le C_\epsilon r^{3/2-\epsilon}.
\end{equation}
These are actual contact estimates, including the vertical bending of the completed cubic.

Suppose in addition that normalized remainders converge at \emph{all} small scales to a fixed two-dimensional slit profile $q(y,z)=c\,\operatorname{Re}(s+i|z|)^{7/2}$, $s=y\cdot e$ for a fixed tangential unit vector $e$, $c>0$.
On the plane the trace on $s>0$ is $c s^{7/2}$; away from the edge, the required normalized convergence is a uniform comparison on compact sets, not just weak convergence.
Along actual contacts on $\partial B_r$, \eqref{access:eq:tiptrace}, \eqref{access:eq:vert}, and the Lipschitz bound on the model imply
\begin{equation}\label{access:eq:ap-from-conv}
 \frac{s_+}{r}\le C_\epsilon
 \left(\|w(r\cdot)/H(r,w)^{1/2}-q\|_{L^\infty(B_2)}
       +r+r^{3/2-\epsilon}\right)^{2/7}.
\end{equation}
Thus mere all-scale uniform convergence, even with the correct unique slit, gives only $s_+/r=o(1)$.
A power rate in this uniform comparison pays a power aperture immediately.
\cite[Proposition 4.12(a)]{FRS} supplies compactness to homogeneous Signorini profiles, not this fixed profile with a rate.
The half-space results of Savin--Yu \cite{SY} are for the thin obstacle problem in three dimensions and do not, without a separate argument, yield \eqref{access:eq:ap-from-conv} with a power rate for this completed remainder of the classical obstacle problem in dimension five.

\subsection{A barrier in the original positivity set}\label{sec:slit-access}

For $S_0=\{(s,z):z=0,\ s\le0\}$ denote by $\rho=(s^2+z^2)^{1/2}$ the two-dimensional distance to its edge.
The nonnegative slit function
\begin{equation*}
 U_{1/2}(s,z)=\operatorname{Re}(s+i|z|)^{1/2}=\sqrt{(\rho+s)/2}
\end{equation*}
is defined with branch cut $S_0$.
It is harmonic off the slit in all five variables, homogeneous of degree $1/2$, and
\begin{equation}\label{access:eq:w}
 |DU_{1/2}|^2=\frac1{4\rho},\qquad x\cdot DU_{1/2}=U_{1/2}/2.
\end{equation}
It satisfies $U_{1/2}^2\le \dist((s,z),S_0)$: for $s\ge0$ use $U_{1/2}^2\le\rho$, and for $s<0$ use $\rho-|s|\le|z|$.

\begin{proposition}[Power aperture implies exact slit access]\label{access:prop:bridge}
Assume for some $\eta>0$ and $C$ the actual contact set satisfies
\begin{equation}\label{access:eq:power}
 \dist((s,z),S_0)\le Cr^{1+\eta}
 \quad\hbox{on }K\hbox{ near the point}.
\end{equation}
For any fixed $\chi>0$, at all sufficiently small radii the actual temporal difference has
\begin{equation*}
 \inf_{\{x\in\partial B_r:\,|z|\ge\chi r\}}
 (u^t(x)-u^0(x))\ge c_\chi t r^{1/2}\qquad (t>0).
\end{equation*}
Consequently the usual signed comparison gives full forward cleaning in $B_{r/2}$ by time $Cr^3$, for this point.
\end{proposition}
\begin{proof}
Let $\alpha=(1+\eta)/2>1/2$, $b=r^\alpha$, $\zeta=U_{1/2}/b$, and choose $C_0>\sqrt C$.
On the open set $\zeta>C_0$ define
\begin{equation*}
 \varphi(\zeta)=s_0-\log(1+s_0),\qquad s_0=\zeta-C_0,
 \qquad Z=b\varphi(\zeta).
\end{equation*}
This set is contained in the actual positivity set by \eqref{access:eq:power}.
The perspective differentiation identity is
\begin{equation}\label{access:eq:perspective}
 \Delta Z=\frac{\varphi''(\zeta)}{b}|DU_{1/2}-\zeta Db|^2
       +(\varphi(\zeta)-\zeta\varphi'(\zeta))\Delta b .
\end{equation}
It follows by differentiating $b\varphi(U_{1/2}/b)$, using $\Delta U_{1/2}=0$.
The first term is nonnegative.
The coefficient in the second is bounded below by $-C_1(1+\log(1+\zeta))$, and $\zeta\le r^{1/2-\alpha}$.
Thus $\Delta Z\ge -C_2r^{\alpha-2}(1+T)$, where $T=\log(1/r)$.
Also
\begin{equation}\label{access:eq:near}
 |DU_{1/2}-\zeta Db|^2=|DU_{1/2}|^2+\alpha(\alpha-1)U_{1/2}^2/r^2,
 \qquad
 \varphi''(\zeta)/b\ge c b/(U_{1/2}+b)^2 .
\end{equation}
Both statements follow from \eqref{access:eq:w} and the explicit formula $\varphi''=1/(1+s_0)^2$.
If $U_{1/2}\le r^{1/2}/T$, then the first expression in \eqref{access:eq:near} is at least $c/r$ when $r$ is small, since $\rho\le r$ and the possible negative term is $O(1/(rT^2))$.
The positive contribution in \eqref{access:eq:perspective} is then at least $c r^{\alpha-2}T^2$.

Choose $0<\beta<\alpha-1/2$ and put $f(r)=1+r^\beta$.
We claim $fZ$ is subharmonic.
The radial derivative is
\begin{equation*}
 \partial_r Z
 =\varphi'(\zeta)\frac{U_{1/2}}{2r}
       +(\varphi(\zeta)-\zeta\varphi'(\zeta))\frac{\alpha b}{r},
\end{equation*}
so the cross term $2f'\partial_r Z$ is bounded below by $-C r^{\alpha+\beta-2}(1+T)$.
On the region $U_{1/2}\le r^{1/2}/T$ the preceding positive contribution dominates this and the negative radial term; $Z\Delta f\ge0$.
On the complementary region, $U_{1/2}/b\gg T$ uniformly as $r\downarrow0$, hence $Z\ge U_{1/2}/2$, and
\begin{equation*}
 Z\Delta f\ge c r^{\beta-3/2}/T .
\end{equation*}
This dominates $r^{\alpha-2}(1+T)$ and the cross-term error since $\beta<\alpha-1/2$.
The claim follows.

The barrier $fZ$ vanishes continuously at $U_{1/2}=C_0b$ and at the tip.
On an outer sphere at a fixed small radius the closure of its support is in the actual warm set, since $C_0>\sqrt C$.
Each point of this compact outer support can be connected within the old positivity set to a transverse cap at a fixed positive radius.
Indeed, along a ray the ratio $U_{1/2}/r^\alpha$ increases as $r$ decreases.
Choose an angular path in $\{U_{1/2}>0\}$ from the ray to a transverse direction; its minimum on the unit sphere is positive, so shrinking first places the entire angular path in $\{U_{1/2}>C_0r^\alpha\}$.
At the cap, \cite[Lemma 9.2(b)]{FRS} gives a lower bound $ct$, with the fixed radius absorbed in $c$.
The increment is nonnegative and harmonic in the old positivity set.
Harnack chains and a finite covering of the compact outer support therefore give the uniform outer bound needed for the comparison.
Compare it with a small multiple of $tfZ$ in the warm comparison domain.
On fixed transverse caps $U_{1/2}\asymp\sqrt r$ uniformly, hence $fZ\asymp\sqrt r$.

Finally take the signed harmonic comparison polynomial from \cite[Lemma 9.1]{FRS}, whose positive caps have $|z|\ge r/\sqrt{10}$ and which is positive in $B_{r/2}$.
The ansatz has $\Delta\mathcal P=1$ and the global lower bound $\mathcal P\ge-Cr^5$ on $B_r$ by \eqref{eq:ansatz-lower}.
The comparison cost $Cr^{7/2}$ comes from \eqref{access:eq:amp}.
With the lower cap increment just proved the time threshold is $C r^{7/2}/r^{1/2}=Cr^3$.
\end{proof}

Proposition~\ref{access:prop:bridge} supplies the analytic cleaning estimate under the contact-aperture hypothesis alone.
Section~\ref{sec:slit-iteration} will obtain that aperture from Theorem~\ref{slit:th:one}, construct the three-dimensional spatial charts separately from the limiting secant directions, and combine the two conclusions to control the parameter projection.

\subsection{Iteration and parameter projection}\label{sec:slit-iteration}
At a fixed base with positive lower density of rich scales for some aperture, choose $p>0$ below that density.
For $r_k=2^{-k}$ there are arbitrarily late starting indices $n$ such that every prefix beginning at $n$ has proportion at least $p$ of rich indices, by the binary good-start lemma \cite[Lemma 5.5]{FS}.
The height ratios $h_0(r_k/2)/h_0(r_k)\to2^{-\lambda}$ follow from normalized homogeneous compactness.
Weaken the strict gain to $\sigma_1<1/64$ and choose the loss $\gamma_1>0$ so small that $p\sigma_1-(1-p)\gamma_1=\kappa>0$.
For sufficiently late scales the normalized recurrence uses factors $2^{-\sigma_1}$ on rich indices and $2^{\gamma_1}$ otherwise, with additive floor $C2^{-k/32}$.
Set $d_k=d(r_k)$ and let $\xi_k=2^{-\sigma_1}$ on rich indices and $\xi_k=2^{\gamma_1}$ otherwise.
For $P_{n,k}=\prod_{i=n}^{k-1}\xi_i$, the good-start property gives
\begin{equation*}
 P_{n,k}\le2^{-\kappa(k-n)},\qquad
 P_{j+1,k}\le P_{n,k}2^{\sigma_1(j+1-n)},
\end{equation*}
where the second inequality uses $\xi_i\ge2^{-\sigma_1}$.
Expanding the recurrence and summing its error terms yields
\begin{align*}
 d_k&\le P_{n,k}d_n+C\sum_{j=n}^{k-1}P_{j+1,k}2^{-j/32}\\
 &\le C2^{-\kappa(k-n)}(d_n+2^{-n/32}),
\end{align*}
because $\sigma_1<1/32$ makes the last geometric series convergent.
The good-start condition for a prefix of length one implies that $n$ itself is rich.
Thus $d_n\to0$ along arbitrarily late good starts by Theorem~\ref{slit:th:one}.
Choose such a start with $C(d_n+2^{-n/32})$ below the smallness threshold.
Induction now justifies every step of the recurrence and keeps all subsequent $d_k$ within that threshold.
Comparing successive optimal profiles on the same smaller ball gives orientation increments bounded by a constant times $d(r_k)+d(r_{k+1})$: use homogeneity of the outer fit and comparable normalizations, then normalize the positive amplitudes.
The pure profile map has locally Lipschitz inverse for orientation on its compact orbit (differentiate in the three tangent directions to see its injective differential).
Thus the slit direction converges at a power rate.
Restriction of the outer fit gives the bound at intermediate scales.
For interior supremum control pass to the corresponding exact thin replacement and use the subharmonic square of its difference with the pure fit.
Estimate \eqref{slit:eq:oddreplacement} bounds the odd part separately.
Thus at every sufficiently small radius, the full normalized remainder differs from a slit of fixed orientation by a positive power of the radius.
The amplitude in this normalization is bounded above and away from zero.
The same comparison yields a nonzero leading coefficient in the original scaling.
If $a_k$ is the positive amplitude of the height-normalized fit, the physical coefficient is $b_k=a_k h_0(r_k)/r_k^\lambda$.
Comparing the unnormalized fits on the same smaller ball shows $|b_{k+1}-b_k|\le C b_k(d(r_k)+d(r_{k+1}))$, using comparable heights.
The resulting summable changes of $\log b_k$ give a positive limit, with a power rate.

For completeness, the spatial projection of the complementary piece has Hausdorff dimension at most two.
For every fixed positive aperture, the two-plane tube property holds with upper density one in dyadic scale at each base of this piece; restricting the set of neighbors to the complementary piece only improves the property.
Apply the purely geometric \cite[Lemma 4.4]{FS}, with the index shift from $B_{r/2}$ immaterial, and send the aperture and density deficits to zero.
The parameter comparison with exponent $5/2-\eta$ \cite[Lemmas 9.1--9.3 and Corollary 7.8(a)]{FRS} gives null parameter projection there.

On the rich piece, the all-scale power comparison just proved and \eqref{access:eq:ap-from-conv} give $s_+/r\le Cr^{\theta}$ for some $\theta>0$.
Together with $|z|\le Cr^2$ from \eqref{access:eq:vert}, this gives
\begin{equation*}
 \dist((s,z),S_0)\le s_++|z|\le Cr^{1+\min(\theta,1)}
 \quad\hbox{on the actual contact set}.
\end{equation*}
Thus the already proved Proposition~\ref{access:prop:bridge} gives $Cr^3$ forward cleaning.
The fixed slit determines a three-plane containing the limits of all spatial matched secants by \cite[Lemma 6.8]{FRS}.
To obtain charts, localize the radius at which every unit secant is within $1/8$ of the limiting three-plane, approximate those planes by a finite $1/8$-net in projection norm, and partition into spatial cubes smaller than the localized radius.
Each piece is a graph over its net plane with Lipschitz inverse projection, because its normalized secants have orthogonal projection to the complement of norm at most $1/4$.
Ordering each pair by time and uniformizing the radii and constants, the projection estimate in Theorem~\ref{sard:th:sard} applies with logarithmic power zero.
This is only the rank-one, harmonic-cubic, next-frequency-$7/2$ cohort.
No nonharmonic cubic or other frequency stratum is included.

\section{Proof of Theorem~\ref{thm:main}}\label{sec:conclusion}
The singular loci where the quadratic kernel has dimension at most two, or dimension three with anomalous second blow-up, have null parameter projection by Section~\ref{sec:reduction}.
The ordinary dimension-three quadratic-kernel locus is covered by Corollary~\ref{sard:cor:closure}.
On the maximal stratum the subcubic interval is empty by the Signorini frequency-gap argument in Section~\ref{sec:projection-reduction}.
If the cubic is nonharmonic, the rank at least two is covered by Section~\ref{sec:reduction}, and rank one by Theorem~\ref{cubic:th:param}.
After a completed harmonic cubic, the next frequency lies in $(3,4)$ or is at least four.
Within $(3,4)$ the frequencies other than $7/2$ are covered by the frequency-localized dimension and cleaning comparison in Section~\ref{sec:reduction}; the exact $7/2$ cohort is handled in Section~\ref{sec:slit}.
At and above four, the failed quartic and completed quartic loci are both treated in Section~\ref{sec:reduction}.
If a frequency normalization fails because the residual is identically zero locally, the corresponding exact nonnegative polynomial has Laplacian one, and the same cap comparison with zero polynomial error gives local forward positivity.
Thus no such point contributes a positive measure of parameters.
Countably exhausting the local cylinders and taking the union proves the theorem.\qed


\begin{thebibliography}{13}

\bibitem{AltafCs}
I. Altaf and M. Cs\H{o}rnyei, On the modulus of continuity of functions
whose image has positive measure, and metric embeddings into
$\mathbb R^d$ without shrinking,
\href{https://arxiv.org/abs/2504.06488}{arXiv:2504.06488} (2025).

\bibitem{CC}
M. Carducci and R. Colombo, Generic regularity of free boundaries
in the obstacle problem for the fractional Laplacian,
\href{https://arxiv.org/abs/2412.16066}{arXiv:2412.16066} (2024).

\bibitem{CV}
M. Carducci and B. Velichkov, An epiperimetric inequality for
odd frequencies in the thin obstacle problem,
\emph{J. Funct. Anal.} \textbf{289} (2025), 111115.
\href{https://doi.org/10.1016/j.jfa.2025.111115}{DOI}; \href{https://arxiv.org/abs/2409.12110v3}{preprint version used for formula numbering}.

\bibitem{FR}
X. Fern\'andez-Real, Existence of odd-frequency free boundary
points in the Signorini problem,
\href{https://arxiv.org/abs/2611.34456}{arXiv:2611.34456} (2026).

\bibitem{FRTL}
X. Fern\'andez-Real and C. Torres-Latorre, Generic regularity of free
boundaries for the thin obstacle problem, \emph{Adv. Math.}
\textbf{433} (2023), 109323.
\href{https://doi.org/10.1016/j.aim.2023.109323}{DOI}.

\bibitem{FRS}
A. Figalli, X. Ros-Oton, and J. Serra, Generic regularity of free boundaries for the obstacle problem,
\emph{Publ. Math. Inst. Hautes \'Etudes Sci.} \textbf{132} (2020), 181--292.
\href{https://doi.org/10.1007/s10240-020-00119-9}{DOI}.

\bibitem{FrSa}
F. Franceschini and O. Savin, A note on the frequency gaps between integers in the thin obstacle problem,
\href{https://arxiv.org/abs/2412.13425}{arXiv:2412.13425} (2024).

\bibitem{FS}
F. Franceschini and J. Serra, Free boundary partial regularity in
the thin obstacle problem, \emph{Comm. Pure Appl. Math.}
\textbf{77} (2024), 630--669.
\href{https://doi.org/10.1002/cpa.22152}{DOI}; \href{https://arxiv.org/abs/2112.11104v1}{preprint version used for lemma numbering}.

\bibitem{FZ}
F. Franceschini and W. Zato\'n, $C^\infty$ partial regularity of the singular set in the obstacle problem,
\emph{Anal. PDE} \textbf{18} (2025), 199--264.
\href{https://doi.org/10.2140/apde.2025.18.199}{DOI}.

\bibitem{Monneau}
R. Monneau, On the number of singularities for the obstacle problem in two dimensions,
\emph{J. Geom. Anal.} \textbf{13} (2003), 359--389.
\href{https://doi.org/10.1007/BF02930701}{DOI}.

\bibitem{SYinteger}
O. Savin and H. Yu, Contact points with integer frequencies in the
thin obstacle problem, \emph{Comm. Pure Appl. Math.} \textbf{76} (2023), 4048--4074.
\href{https://doi.org/10.1002/cpa.22126}{DOI};
\href{https://arxiv.org/abs/2103.04013}{arXiv:2103.04013}.

\bibitem{SY}
O. Savin and H. Yu, Half-space solutions with $7/2$ frequency in the thin obstacle problem, \emph{Arch. Ration. Mech. Anal.} \textbf{246} (2022), 397--474.
\href{https://doi.org/10.1007/s00205-022-01817-w}{DOI:10.1007/s00205-022-01817-w}.

\bibitem{Schaeffer}
D. G. Schaeffer, An example of generic regularity for a nonlinear elliptic equation,
\emph{Arch. Ration. Mech. Anal.} \textbf{57} (1974), 134--141.
\href{https://doi.org/10.1007/BF00248415}{DOI}.

\end{thebibliography}
\end{document}